\documentclass[10pt]{article}
\usepackage{amscd,amsfonts,amsmath,amssymb,amstext,amsthm}
\newcommand{\forgetIt} [1] { } 
\newcommand{\wrt} {with respect to\ } 
\newcommand{\s}{\par\smallskip \noindent}
\newcommand{\m}{\par\medskip \noindent}
\renewcommand{\b}{\par\bigskip \noindent}
\theoremstyle{plain}
\newtheorem{thm} {Theorem} [subsection]
\newtheorem{lem} [thm]{Lemma}
\newtheorem{prop}[thm]{Proposition}
\newtheorem{cor} [thm]{Corollary}
\theoremstyle{definition}
\newtheorem*{defn} {Definition}

\theoremstyle{remark}
\newtheorem*{rem} {Remark}

\font    \sixrm=cmr6      
\font    \fourrm=cmr5  scaled 800
\newcommand{\OAG} {{\Omega_{\kern-1pt\rm AG}}}
\newcommand{\oag} {{\omega_{\kern-1pt\rm AG}}}
\newcommand{\Aut} {\mathop{\rm Aut}\nolimits}
\newcommand{\GL}  {\mathop{\rm GL}\nolimits}
\renewcommand{\P}{{\rm I\kern-.19em P}}
\def\sL{{\mathfrak s}{\mathfrak l}}
  
\def\ss{{\rm ss}}
\def\buildre#1\over#2{\mathrel{\mathop{\kern0pt #1}\limits_{#2}}}
\def\GCKC{G^{\C}\kern-1.5pt /\kern-1pt K^{\C}}
\newcommand{\Ad}{\mathop{\rm Ad}\nolimits}
\newcommand{\codim}{\mathop{\rm codim}\nolimits}
\newcommand{\Int}{\mathop{\rm Int}\nolimits}
\newcommand{\diag}{\mathop{\rm diag}\nolimits}
\newcommand{\End}{\mathop{\rm End}\nolimits}
\newcommand{\Cl}{{\cal C}\kern-1pt\ell}
\def\k{{\mathfrak k}}
\def\h{{\mathfrak h}}
\def\g{{\mathfrak g}}

\def\l{{\mathfrak l}}
\def\p{{\mathfrak p}}
\def\q{{\mathfrak q}}

\def\ss{{\mathfrak s}}
\def\mm{{\mathfrak m}}

\def\gc{{\mathfrak g}^{\C}}

\def\kc{{\mathfrak k}^{\C}}
\def\ac{{\mathfrak a}^{\C}}

\def\sc{{\mathfrak s}^{\C}}

\def\pc{{\mathfrak p^{\C}}}
\def\R{\mathbb R}
\def\Z{\mathbb Z}
\def\C{\mathbb C}
\def\a{{\mathfrak a}}

\def\OA{\Omega_{\kern-1pt \rm AG}}\def\oa{\omega_{\kern-1pt\rm AG}}

\def\gc{{\mathfrak g}^\C}\def\kc{{\mathfrak k}^\C}\def\pc{{\mathfrak p}^\C}
 \def\sc{{\mathfrak s}^\C}

\catcode`\@=11
\def\diagram#1{\lineskiplimit=0pt\null\,\vcenter{\openup\jot\m@th
  \ialign{&\strut$\displaystyle\hfil{}##{}\hfil$\crcr#1\crcr}}\,}%
\catcode`\@=12
\def\R{{\mathchoice
{\hbox{\rm I\kern-.185em R}} {\hbox{\rm I\kern-.185em R}}
{\hbox{\sixrm I\kern-.19em R}}{\hbox{\fourrm I\kern-.18em R}}
}}
\def\C{{\mathchoice{
\hbox{\rm\kern.32em\vrule height1.5ex width.07em depth-.011em\kern-.34em C}}
{\hbox{\rm\kern.32em\vrule height1.5ex width.07em depth-.011em\kern-.34em C}}
{\kern-.12em\hbox{\sixrm\kern.39em\vrule height1.48ex width.065em depth-.023em\kern-.35em C}\kern-1pt}
{\kern-.13em\hbox{\fourrm\kern.37em\vrule height1.47ex width.06em depth-.023em\kern-.32em C }\kern-3pt}}}
\title {Characterization of cycle domains via\\
Kobayashi hyperbolicity}
\author {Gregor Fels\footnote{Research supported by a 
Habilitationsstipendium from the Deutsche Forschungsgemeinschaft}\ \  and 
Alan Huckleberry\footnote {Research partially supported by grants
from the Deutsche Forschungsgemeinschaft and the Japanese Society for the
Promotion of Science}}
\date{}
\begin{document}
\maketitle
\abstract{\s
\begin{quote}\footnotesize
A real form $G$ of a complex semisimple Lie group $G^\mathbb C$ has only
finitely many orbits in any given $G^\mathbb C$-flag manifold
$Z=G^\mathbb C/Q$.  The complex geometry of these orbits is of interest, e.g.,
for the associated representation theory.  The open orbits $D$ generally
possess only the constant holomorphic functions, and the relevant associated
geometric objects are certain positive-dimensional compact complex 
submanifolds of $D$ which, with very few well-understood exceptions,
are parameterized by the Wolf cycle domains
$\Omega _W(D)$ in $G^\mathbb C/K^\mathbb C$, where $K$ is a maximal
compact subgroup of $G$.  Thus, for the various domains $D$ in
the various ambient spaces $Z$, it is possible to compare the cycle 
spaces $\Omega _W(D)$.  
\s
The main result here is that, with the few exceptions mentioned above, for 
a fixed real form $G$ all of the cycle spaces $\Omega _W(D)$ are
the same.  They are equal to a universal domain $\Omega _{AG}$ which
is natural from the the point of view of group actions and which, 
in essence, can be explicitly computed.
\s
The essential technical result is that if $\widehat \Omega $ is a
$G$-invariant Stein domain which contains $\Omega _{AG}$ and which
is Kobayashi hyperbolic, then $\widehat \Omega =\Omega _{AG}$. The
equality of the cycle domains follows from the fact that 
every $\Omega _W(D)$ is itself Stein, is hyperbolic, 
and contains $\Omega _{AG}$.\end{quote}}
%
%

\section {Introduction} \label {I}

Let $G$ be a non-compact real semi-simple Lie group which is
embedded in its complexification $G^{\C}$ and consider the
associated $G$-action on a $G^{\C}$-flag manifold
$Z=G^{\C}/Q$.  It is known that $G$ has only finitely
many orbits in $Z$; in particular, there exit open $G$-orbits
$D$.  In each such open orbit every maximal compact subgroup $K$ of 
$G$ has exactly one orbit $C_0$ which is a complex submanifold
(\cite {W1}).
\m
Let $q:=\dim_{\C} C_0$, regard $C_0$ as a point in the space
${\cal C}^q(Z)$ of $q$-dimensional compact cycles in $Z$ and let
$\Omega :=G^{\C}\cdot C_0$ be the orbit in ${\cal C}^q(Z)$. 
Define the Wolf cycle space $\Omega _W(D)$ to be the connected component of
$\Omega \cap {\cal C}^q(D)$ which contains the base cycle
$C_0$.
\m
Since the above mentioned basic paper (\cite {W1}) there has
been a great deal of work aimed at describing 
these cycle spaces.  Even in  situations where good
matrix models are available this is not a simple matter.
In fact an exact description of $\Omega _W(D)$ 
has only been given in very special situations 
(see e.g. \cite {BLZ}, \cite {BGW}, \cite {DZ}, \cite {HS},
\cite {HW1}, \cite {N}, \cite {PR}, \cite {We}, \cite {WZ2}). 
In concrete cases of complex geometric relevance, such as the 
period domain $D$ for marked $K3$-surfaces, only partial information was
available (see \cite {Hr},\cite {H}).  Recent progress,
which is outlined below, only yielded qualitative information,
such as the holomorphic convexity of $\Omega _W(D)$
(\cite {W2},\cite {HW2}).  
\m
The results in the present paper change 
the situation: Except for a few exceptional cases, which we discuss below,
for a given group $G$ the Wolf cycle spaces $\Omega_W(D)$
are all the same and are equal to a domain $\OAG$ 
which can be explicitly described.  The main new ingredients here
involve a combination of complex geometric and combinatorial
techniques for studying $G$-invariant, Kobayashi hyperbolic, 
Stein domains in the associated cycle space $\Omega $.
\s
As indicated above, we shall now outline a number of recent developments
and then state our main results. In order to avoid cumbersome 
statements involving products, it will be assumed that $G$ is simple.  
This is no loss of generality.
\m
In (\cite {W2}) Wolf shows that
if the orbit $\Omega=G^{\C}\!\cdot C_0$ in the cycle space is 
{\it compact\/}, then it is either of Hermitian type or just a point.
In the former case $G$ is of Hermitian type and
the cycle space $\Omega_W(D)$ is the dual bounded symmetric
domain ${\cal B}\subset\Omega$ (or its complex conjugate
$\bar {\cal B}$). The later case occurs if the real form 
$G_0$ acts transitively on
$Z$ (see \cite {W3}, \cite {On} for a classification).
\s
Here we do not consider the above described 
well-understood special cases but investigate only the generic 
case, where $\Omega$ is  {\it non-compact}. In this case it is 
a complex-affine variety such that the connected component of 
the isotropy group $G^{\C}_{C_0}$ is exactly the 
complexification $K^{\C}$.
\s
We replace $\Omega =G^\C/G^\C_{C_0}$ by the finite
covering space $G^{\C}/K^{\C}$ (which in many cases is the same as
$\Omega $), choose a base point $x_0$ with isotropy $K^{\C}$ and regard 
$\Omega _W(D)$ as the "orbit of $x_0$" of the connected component of
$\{ g\in G^{\C}:g(C_0)\subset D\}$ which contains the identity.
In this way it is possible to compare all cycle domains for a
fixed group $G$.
\b
One of the main motivating factors for studying the complex geometry
of the cycle domains $\Omega _W(D)$ is that ${\cal O}(\Omega _W(D))$,
or, more generally, spaces of sections of holomorphic vector bundles
over $\Omega _W(D)$, provides a rich source of $G$-representations.
\s
As a general principle one relates representation theory and complex analysis
by starting with a smooth manifold $M$ equipped, e.g., with a proper
$G$-action and $G$-equivariantly complexify $M$ to obtain a holomorphic
$G^{\C}$-action on a Stein manifold $M^\C$ and a 
$G$-invariant basis of Stein neighborhoods of $M$ where the $G$-action
is still proper, properness being useful for constructing invariant
metrics, volume forms, etc..  This type of complexification exists 
in general (see \cite {He},\cite {Ku},\cite {HHK}), but a concrete description of, e.g., 
{\it maximal\/} $G$-invariant Stein neighborhoods where the $G$-action
remains proper is a difficult matter in the general setting.
\m
Here we consider the special case of $M=G/K$ a Riemannian symmetric
space of non-positive curvature (and $G$ semisimple) embedded in 
$M^{\C}=G^{\C}/K^{\C}$ as an orbit of the same
base point $x_0$ as was chosen above in the discussion of cycle spaces.
In this setting Akhiezer and Gindikin (\cite {AG}) define a 
neighborhood $\OAG$ of $M$ in $M^{\C}$ which is quite
natural from the point of view of proper actions.
\m
For this we fix a bit of the notation which will be used throughout
this paper.  Let $\theta $ be a Cartan involution on 
${\mathfrak g}^{\C}$ which restricts to a Cartan involution on
${\mathfrak g}$ such that $Fix(\theta \vert {\mathfrak g})={\mathfrak k}$
is the Lie algebra of the given maximal compact subgroup $K$.
The anti-holomorphic involution 
$\sigma :{\mathfrak g}^{\C}\to {\mathfrak g}^{\C}$ 
which defines ${\mathfrak g}$
commutes with $\theta $ as well as with the holomorphic extension
$\tau $ of $\theta \vert {\mathfrak g}$ to ${\mathfrak g}^{\C}$.
\m
Let ${\mathfrak u}$ be the fixed point set of $\theta $ in
${\mathfrak g}^{\C}$, $U$ be the associated maximal
compact subgroup of $G^{\C}$ and define $\Sigma $ 
to be the connected component containing $x_0$ of 
$\{ x\in U\cdot x_0: G_x\ \text{\rm is compact}\} $.
Set $\OAG:=G\cdot \Sigma $.
It is shown in (\cite {AG}) that $\OAG$ is an open
neighborhood of $M:=G\cdot x_0=G/K$ in $\GCKC$ on which the action of $G$ is
proper.  
\m
To cut down on the size of $\Sigma $, one considers a maximal Abelian  
subalgebra ${\mathfrak a}$ in 
${\mathfrak p}$ (where $\g=\k\oplus \p$ is the Cartan decomposition of $\g$) 
and notes that 
$G\cdot (\exp(i{\mathfrak a})\cap \Sigma )\cdot x_0=\OAG$.
In fact there is an  explicitly defined neighborhood 
$\oag$ of $0\in {\mathfrak a}$ such that 
$i{\mathfrak \oag }$ is mapped diffeomorphically onto its
images $\exp (i\oag)$ and $\exp(i\oag )\cdot x_0$  and 
$\OAG=G\cdot \exp(i\oag)\cdot x_0$.  
The set $\oag $ is defined by the 
set of roots $\Phi ({\mathfrak a})$ of the
adjoint representation of ${\mathfrak a}$ on ${\mathfrak g}$:
It is the connected component containing $0\in {\mathfrak a}$
of the set which is obtained from $\a$ by removing the root hyperplanes
$\{ \mu=\pi/2 \} $ for all $\mu\in \Phi(\a) $. It is convex and
is invariant under the action of the Weyl group ${\cal W}{(\mathfrak a})$ 
of the symmetric space $G/K$. Modulo 
${\cal W}({\mathfrak a})$, the set $\exp(i\oag)\cdot x_0$ is
an geometric slice for the $G$-action on $\OAG$;
(see \cite {AG} for details).
\b
A further viewpoint arises in the study of Schubert incidence varieties
in which  $\Omega =G^{\C}/K^{\C}$ is again regarded as a space of 
cycles (\cite {HS},\cite {H}, \cite {HW1}, \cite {HW2}).  For this,
a Borel subgroup $B$ of $G^{\C}$ which contains the factor
$AN$ of an Iwasawa-decomposition $G=KAN$ is called a Borel-Iwasawa
group and a $B$-Schubert variety $S\subset Z=G^{\C}/Q$ an
Iwasawa-Schubert variety, i.e., the closure of some orbit of a
Borel-Iwasawa subgroup.
\m
Given a boundary point $p\in \text{\rm bd}(D)$, the incidence variety
approach which was exemplified in (\cite {HS}) and developed to
its completion in (\cite {H}) and (\cite {HW2}), yields an
Iwasawa-Schubert variety $Y$ with $\codim_{\C}Y=q+1$ contained in
the complement $Z\smallsetminus D$ with $p\in Y$.  The incidence
hypersurface
$H_Y:=\{ C\in \Omega :C\cap Y\not=\emptyset \} $ lies in the complement of
$\Omega _W(D)$ and in particular touches its boundary at the cycles which
contain $p$. If $H=H_Y$ is moved by $g\in G$, then of course 
$g(H)\cap \Omega _W(D)=\emptyset $ as well.  
\s
Thus, for any
hypersurface $H$ in $\Omega $ which is invariant by some Borel-Iwasawa subgroup
$B$ it is appropriate to consider the domain 
$\Omega _H$ which is the connected component containing
the base point $x_0$ of $\Omega \smallsetminus {\cal E}(H)$, where ${\cal E}(H)$ 
is the closure of the union of the hypersurfaces $g(H)$, $g\in G$.
\m
Since $\Omega $ is a spherical homogeneous space, a Borel subgroup
$B$ has only finitely many orbits in $\Omega $.  In particular,
there are only finitely many $B$ invariant hypersurfaces
$H_1,\ldots ,H_m$.  For a given Iwasawa-Borel subgroup let
$I(D):=\{ j:H_j\subset \Omega \smallsetminus \Omega _W(D)\}$ and 
$I:=\{ 1,\ldots ,m\} $ the full index set.
\m
The Schubert domain $\Omega _S(D)$ is defined to be the connected
component containing $\Omega _W(D)$ of the intersection
$\bigcap _{j\in I(D)}\Omega _{H_j}$ and the universal Iwasawa domain
$\Omega _I$ to be the connected component containing the base point
$x_0$ of $\bigcap _{j\in I}\Omega _{H_j}$.  Of course all of these
domains are $G$-invariant and Stein.  The above mentioned existence of 
Iwasawa-Schubert varieties $Y$ at boundary points of $D$ implies that
$\Omega _W(D)=\Omega _S(D)$, one of the main results in (\cite {HW2}). 
On the analytic side the domain $\Omega _I$ plays a role in 
matters of holomorphic extension of certain spherical functions 
on $G/K$ (see \cite{Fa}, \cite{Ola}, \cite{OO}, \cite{KS}).
\m
It is known that 
$\Omega _I=\OAG$ (see $\S$\ref {II}). 
Hence, excluding the few exceptional cases mentioned above,
before the present paper it was known that
$\OAG=\Omega _I\subset \Omega _S(D)=\Omega _W(D)$ 
for all open $G$-orbits $D$. It should be underlined
that apriori, for a fixed real form $G$ but arbitrary open
$G$-orbits $D$ in the various projective $G^\C$-homogeneous manifolds  
$Z$, the cycle spaces $\Omega _W(D)$ might all be different.
\m
Here, with the few well-understood exceptions mentioned above, we
prove the following equality.
\m
{\bf Theorem}(see \ref{description}): {\it $\OAG=\Omega _S(D)$ for any open 
$G$-orbit $D$ in
any $G^{\C}$-flag manifold. }
\m
Thus, $\Omega _W(D)$ is either a point, a Hermitian bounded symmetric space
or (the generic case) $\Omega _W(D)=\OAG $.
\s
Our main results are consequences of considerations involving the Kobayashi
pseudo-metric.  In particular, building on work from 
(\cite {H}), we show that the domains $\Omega _S(D)$ 
are Kobayashi hyperbolic (see $\S$\ref {V}).  It is clear from 
the definitions that they are $G$-invariant and Stein. The essential
point can then be formulated as follows 
\m
{\bf Theorem}(see \ref {notbrody}):
{\it If $\widehat \Omega $ is a $G$-invariant, Kobayashi hyperbolic,
Stein domain containing $\OAG$, then $\widehat \Omega =\OAG$.}
\m
The equality  $\OAG=\Omega _S(D)=\Omega _W(D)$ is a 
direct corollary of this theorem.


\section {Equivalence of $\OAG$ and $\Omega _I$} \label {II}

Here we regard $\Omega _I$ and $\OAG$ as being set-up at
the chosen base point $x_0$ in $G^{\C}/K^{\C}$ and
show that they are in fact equal.  As one would expect, this
involves proving two inclusions. Partial results have been obtained by
various authors from various points of view (see \cite {B}
and \cite {KS}). In the sequel we give the
complex analytic proof of the inclusion $\OAG\subset \Omega _I$
from (\cite {H}) and an elementary proof of the other inclusion which
in certain aspects follows that in (\cite {B}). 
\m
For the inclusion $\OAG\subset \Omega _I$ the following
observation on the critical set of plurisubharmonic function
is of use.
 
\begin {lem}
 
Let $\rho $ be a strictly plurisubharmonic function on a
$n$-dimensional complex manifold $X$ and $M$ be a real
submanifold of $X$ with $M\subset \{ d\rho =0\}$. Then $M$
is isotropic with respect to the symplectic form
$dd^c\rho $; in particular it is at most $n$-dimensional and
is totally real.
 
\end {lem}
 
\begin {proof}
 
If $\lambda :=d^c\rho $ and $i:M\hookrightarrow X$ is the natural 
injection, then $ \imath ^*(\lambda )=0$.  Therefore
$i^*(dd^c\rho )=di^*(\lambda )=0$.   
\end{proof}
\s
The inclusion $\OAG\subset \Omega _I$ follows
from the existence of a $G$-invariant strictly plurisubharmonic
function on $\OAG$.
           
\begin {prop}

If $\OAG$ and $\Omega _I$ are set-up at the base
point $x_0$, then
$$\OAG\subset \Omega _I.$$

\end {prop}

\begin {proof}
 
Let $\Omega _{adpt}$ be the maximal domain of existence of the
adapted complex structure in the tangent bundle
$G\times _K{\mathfrak p}$ of the Riemannian symmetric space
$G/K$.  It is known that the polar map
$\Omega _{adpt}\to G^{\C}/K^{\C}$, $[(g,\xi)]\mapsto
g\exp(i\xi )$, is a biholomorphic map onto $\OAG$ (\cite {BHH}).
\s
The square of the norm function on the tangent bundle of $G/K$ which is
defined by the invariant metric is strictly plurisubharmonic in
the adapted structure.  Transporting it via the polar map,
we have a $G$--invariant strictly plurisubharmonic function
$\rho :\OAG\to {\R}^{\ge 0}$. Note that $\{ \rho = 0\}=\{ d\rho =0\} $
is the totally real orbit $G\cdot x_0$.
\m        
Now let $B$ be an Borel-Iwasawa subgroup of $G^{\C}$.
If $\Omega _I$ would not contain $\OAG$, there would
exit a $B$-invariant hypersurface $H$ in $G^{\C}/K^{\C}$
having non-empty intersection with $\OAG$. This cannot happen for the following
reason:
\m
Recall that $B$ contains an Iwasawa-component $AN.$ 
Observe also that, since $\OAG$ can
be identified with a domain in the tangent bundle of $G/K$, the action
of $AN$ on $\OAG$ is free and the $AN$--orbits are parameterized by $\Sigma:=K\exp(i\oag)\cdot x_0$.
In particular, $\dim_{\R} AN\cdot y=\dim_{\C} \OAG$.
Therefore, since $\dim_{\C}(H\cap \OAG)= \dim_{\C} \OAG-1$, the $AN$-orbits in $H\cap \OAG$ are not totally
real.
\s
Now, on the other hand $AN\cdot x_0=G\cdot x_0\cong G/K$ is 
a totally real submanifold of $\OAG $. Therefore, for
$x\in exp(i\oag ).x_0$ near $x_0$, so is the orbit $AN\cdot x$. 
Consequently, if $\rho (x)$ is sufficiently small, then
$AN\cdot x$ is still totally real.
Let $r_0$ be the smallest value of $\rho (y)$, for $y\in H\cap \Sigma $.
It follows that $\rho (y)>0$. Now, the restriction
of $\rho $ to the hypersurface  $X:=A^\C N^\C\cdot y \cap \OAG $ has
a minimum at least along $M=AN\cdot y$.
But this is contrary to the above Lemma, because 
$AN\cdot y$ is not totally real.
\end{proof}
\s
The inclusion $\Omega_I\subset \OAG$ follows directly from our main
theorem \ref{description}. Here, we give a another completely elementary
proof. 

\begin {thm}  \label {equality}

$\OAG=\Omega _I.$

\end {thm}

\begin {proof} 

It is only necessary to prove the opposite inclusion to the 
one above. For this, assume that $\Omega _I$ is not contained
in $\OAG$.  Then there exits a sequence 
$\{ z_n\}\subset \OAG\cap \Omega _I$ with
$z_n\to z\in \text{\rm bd}(\OAG)\cap \Omega _I$.
From the definition of $\OAG$, it follows that
there exist $\{ g_m\}\subset G$ and 
$\{ w_m\} \subset \exp(i\oag )$ such that $g_m(w_m)=z_m$.
\m
Write $g_m=k_ma_mn_m$ in a $KAN$-decomposition of $G$. Since
$\{ k_m\} $ is contained in the compact group $K$, it may be
assumed that $k_m\to k$; therefore that $g_m=a_mn_m$.
\s
Since $\oag $ is relatively compact in ${\mathfrak a}$, it
may also be assumed that $w_m\to w\in \text{\rm c}\ell(\exp(i\oag ))$.
Thus $w_m=A_mx_0$, where $\{ A_m\} \subset \exp(i\oag )$
and $A_m\to A$.  Write $a_mn_n(w_m)=a_mn_mA_mx_0=
\tilde a_m\tilde n_mx_0$, where $\tilde a_m=a_mA_m$
and $\tilde n_m=A_m^{-1}n_mA_m$ are elements of $A^{\C}$ and
$N^{\C}$, respectively.
\m
Now $\{ z_m\}$ and the limit $z$ are contained in $\Omega _I$
which is in turn contained in $A^{\C}N^{\C}\cdot x_0$.
Furthermore, $A^{\C}N^{\C}$ acts freely on this orbit.
Thus $\tilde a_m\to \tilde a\in A^{\C}$ and
$\tilde n_m\to \tilde n\in N^{\C}$ with 
$\tilde a\tilde n\cdot x_0=z$.  Since $A_m\to A$, it follows that
$a_m\to a\in A$ and $n_m\to n\in N$ with $an.w=z$.
Since $z\notin \OAG$, it follows that 
$w\in \text{\rm bd}(\exp(i\oag ))$, and $z\in \Omega _I$
implies that $w\in \Omega _I$.
\m
On the other hand, since $w\in bd(\exp(i\oag ))$, the
isotropy group $G_w$ is non-compact.  But $\Omega _I$ is
Kobayashi hyperbolic (\cite {H}).  Therefore
the $G$-action on $\Omega _I$ is proper (see e.g. \cite {H}) and 
consequently $w\notin \Omega _I$, which is a contradiction.
\end{proof}
\s
We note that the standard definition of $\Omega _I$
looks somewhat different from the one above (see e.g. \cite {B}). 
However, the two definitions are, in fact, equivalent (see \cite {HW2}).


\section {Spectral properties of $\OAG$} \label {III}

\subsection {Linearization}  \label {III.1}

The map $\eta :G^{\C}\to 
\Aut_{\R}({\mathfrak g}^{\C})$, 
$x\mapsto \sigma \circ \Ad(x)\circ \tau \circ \Ad(x^{-1})$,
provides a suitable linearization of the setting at hand
(see \cite {M2} for other applications of $\eta $). In this
section  basic properties of $\eta $ are summarized.
\m
Let $G^{\C}$ act on $\Aut_{\R}({\mathfrak g}^{\C})$
by $h.\varphi :=\Ad(h)\circ \varphi \circ \Ad(h^{-1})$.

\begin {lem}

{\bf ($G$-equivariance)} For $h\in G$ 
it follows that $\eta (h\cdot x)=h. \eta (x)$ for all $x\in G^{\C}$.

\end {lem}

\begin{proof}

By definition $\eta (h\cdot x)=\sigma \Ad(h)\Ad(x)\tau \Ad(x^{-1})\Ad(h^{-1})$.
Since $h\in G$, it follows that $\sigma $ and $\Ad(h)$ commute and
the desired result is immediate.
\end{proof}
\s
The normalizer of $K^{\C}$ in $G^{\C}$ is denoted by
$N^{\C}:=N_{G^{\C}}(K^{\C})$. It is indeed the
complexification of $N:=N_U(K)$.

\begin {lem}

{\bf ($N^{\C}$-invariance)} The map $\eta $ factors
through a $G$-equivariant embedding of 
$G^{\C}/N^{\C}$:

$$\eta (x)=\eta (y) \iff y=xg^{-1}\ \text{\rm for some
}g\in N^{\C}.$$

\end {lem}

\begin{proof}

We may write $y=xg^{-1}$ for some $g\in G^{\C}$.  Thus it must
be shown that $\eta (x)=\eta (xg^{-1})$ if and only if $g\in N^{\C}$.
But $\eta (x)=\eta (xg^{-1})$ is equivalent to $\Ad(g)\tau =\tau \Ad(g),$
which, in turn, is equivalent to the fact that $\Ad(g)$ stabilizes
the complexified Cartan decomposition 
$\gc=(\gc)^\tau\oplus (\gc)^{-\tau}=\kc\oplus \pc$.
\s
Now, if $\Ad(g)$ stabilizes $\kc\oplus \pc$, then 
$\Ad(g)(\kc)=\kc$, i.e., $g\in  N^{\C}.$
On the other hand, given any $g\in  N^{\C}$, it follows $\Ad(g)(\pc)=\pc$,
because $\pc$ is the orthogonal complement of 
$\kc$ \wrt the Killing form of $\gc$.
\end{proof}
\s
Note that $N^{\C}/K^{\C}$ is a finite Abelian group 
(see \cite {Fe} for a classification).
Consequently, up to finite covers, $\eta $ is an embedding of the basic space
$G^{\C}/K^{\C}$.
\m
The involutions $\sigma $ and $\tau $ are regarded as acting
on $\Aut_{\R}({\mathfrak g}^{\C})$ by conjugation.
On $Im(\eta )$ their behavior is particularly simple.

\begin {lem} \label {involutions}

{\bf (Action of the basic involutions)} For all $x\in G^{\C}$
it follows that

\begin{enumerate}

\item $\eta (\tau (x))=\tau(\eta (x))$

\item $\sigma (\eta (x))=\eta (x)^{-1}$

\end{enumerate} 
\s
In particular $Im(\eta )$ is both $\sigma $- and $\tau $-invariant.

\end {lem}

\begin{proof}

Let $\varphi_*:\gc\to \gc$ denote the differential of 
$\varphi:G^\C\to G^\C$ and
$\Int(x):G^\C\to G^\C$ be defined by $\Int(x)(z):=xzx^{-1}.$
The first statement follows directly from the facts that  
$\sigma $ and $\tau $ commute and
$$\tau \Ad(x)\tau =(\tau \Int(x)\tau )_*=\Int(\tau (x))_*=\Ad(\tau (x)).$$
For the second statement note that $\sigma \eta (x)=\Ad(x)\tau \Ad(x^{-1})$,
and thus $\eta (x)\sigma \eta (x)=\sigma $.
\end{proof}
\s
We have seen that $\eta $ is a $G$-equivariant map which induces a 
finite equivariant map 
$\eta:\GCKC\to \Aut_\R(\gc).$ We will shortly see that the image 
$\eta(\GCKC)$ is also closed in $\Aut_\R(\gc).$ Hence, for a characterization 
of $G$--orbits in $\GCKC$ and their topological properties 
we may identify $\GCKC$ with its image in $\Aut_\R(\gc)$ on which 
$G$ acts by conjugation.
\s
The following special 
case of a general result on conjugacy classes (see \cite {Hu} p.~117 and \cite {Bir})
is of basic use.

\begin {lem} \label {ss then closed}

Let $V$ be a finite-dimensional ${\R}$-vector space,
$H$ a closed reductive algebraic subgroup of $\GL _{\R}(V)$ and
$s\in \GL _{\R}(V)$ an element which normalizes $H$. Regard
$H$ as acting on $\GL _{\R}(V)$ by conjugation. Then, 
for a semisimple  $s$  the orbit $H.s $ is closed. 

\end {lem}

\begin {cor} 

The image $Im(\eta )$ is closed in
$\Aut_{\R}({\mathfrak g}^{\C})$.

\end {cor}

\begin {proof}

It is enough to show that $G^{\C} .\tau 
=\{ \Ad(g)\tau \Ad(g^{-1}):g\in G^{\C}\} $ is closed
in $\Aut_{\C}({\mathfrak g}^{\C})$.  Since
$\tau $ is semi-simple and normalizes $G^{\C}$ in
this representation, this follows from Lemma \ref {ss then closed}.
\end{proof}
\s

\subsection {Jordan Decomposition} \label {III.2}

Here $x$ denotes an arbitrary element of $G^{\C}$ and
$su=us=\eta (x)$ is its Jordan decomposition in 
$\GL_{\R}(\gc)$.  Since $\Aut_\R(\gc)$ is algebraic,
$s,u\in \Aut_{\R}({\mathfrak g}^{\C})$ as well.
If $\eta (x)=us$ is not semi-simple, i.e., $u\ne 1$,  consider $\xi =\log(u)\in\End_\R(\gc).$
Since $\xi$ is nilpotent, $t\mapsto \exp(t\xi )$ is an algebraic map
and $\exp(\Z\xi)\subset \Aut_{\R}(\gc).$ 
It follows that
$\exp(t\xi )\in \Aut_{\R}({\mathfrak g}^{\C})$ for all
$t\in {\R}$. In particular, $u$ is in the connected component
$\Aut_{\R}({\mathfrak g}^{\C})^\circ $, and 
 $\xi$ is a derivation: $\xi={\rm ad}(N)$ for some nilpotent $N\in \gc.$ Finally,
 $u=\Ad(\exp(N))=\exp({\rm ad}\, N)$.
\s
Given an element $z\in \Aut_\R(\gc),$ let $(\gc)^z=\{X\in \gc : z(X)=X\}$ denote the subalgebra of fixed points. Observe also that if $\nu:\gc\to \gc$ is any involution such that
$\nu(z)=z$ or $\nu(z)=z^{-1}$, then the subalgebra $(\gc)^z$ is $\nu$-stable.
For $z$ semisimple the subalgebra $(\gc)^z $ is reductive. 

\begin {prop}   \label {lifting}

{\bf (Lifting of the Jordan decomposition)} For $x\in G^{\C}$ with
Jordan decomposition $\eta (x)=u\cdot s$ there exists a nilpotent element
$N\in ({\mathfrak g}^{\C})^s\cap i{\mathfrak g}$
such that

\begin {enumerate}

\item $u=\Ad(\exp(N))$

\item $\eta (\exp(\frac{1}{2}N)\cdot x)=s$.

\end {enumerate}

\end {prop}

\begin{proof}

Let $N\in \gc$ be the element with $u=\Ad(\exp(N))$ as explained above.
First we show that $N\in ({\mathfrak g}^{\C})^s\cap i{\mathfrak g}.$
From Lemma \ref {involutions}(2) it follows that $\sigma (\eta (x))=\sigma(us)=s^{-1}u^{-1}.$
This implies
$\sigma (u)=u^{-1}$ or, equivalently, $\sigma (N)=-N$, i.e.,
$N\in i{\mathfrak g}$. Secondly, the statement that 
$\Ad(\exp(N))$ commutes with $s$ is equivalent to 
$se^{ad(N)}s^{-1}=e^{ad(N)}$ which is the same as $s(N)=N$
in the semi-simple case.  Thus 
$N\in ({\mathfrak g}^{\C})^s\cap i{\mathfrak g}$.
\s
Finally, since $N\in (\gc)^s$, it follows that $\Ad(\exp({t}N))$
commutes with $s$ for all $t\in \R$.
Having also in mind that
$\sigma (N)=-N$, it follows that
 \begin{eqnarray*}
 \eta(\exp \textstyle{\frac{1} {2}} N\cdot x)&=&\sigma\Ad(\exp \textstyle{\frac{1} {2}} N)\Ad(x)\tau\Ad(x^{-1})\Ad(\exp -\textstyle{\frac{1} {2}}N)\vphantom{\big|} = \cr
&=& \Ad(\exp -\textstyle{\frac{1} {2}}N)\cdot \sigma\Ad(x)\tau\Ad(x^{-1})\cdot \Ad(\exp -\textstyle{\frac{1} {2}}N)\vphantom{\raise2pt\hbox{\Big\vert}} =\cr
&=& \Ad(\exp -\textstyle{\frac{1} {2}}N)\cdot s u \cdot \Ad(\exp -\textstyle{\frac{1} {2}}N)=\vphantom{\Big|} \cr
&=& s \cdot \Ad(\exp -\textstyle{\frac{1} {2}}N)\cdot u \cdot 
\Ad(\exp -\textstyle{\frac{1} {2}}N)=s.\vphantom{\big|}
\end{eqnarray*}
\end{proof}
\s
Observe now that since $\sigma(s)=s^{-1},$ $(\gc)^s$ is a $\sigma$--stable 
reductive subalgebra. Let $(\gc )^s=
{\mathfrak h}\oplus {\mathfrak q}$ be its $\sigma $-eigenspace
decomposition.
We now build 
an appropriate ${\mathfrak s\mathfrak l}_2$-triple  $(E,H,F)$ around $N=E$ in
$({\mathfrak g}^{\C})^{s}$.

\begin {lem} \label {sltriple}
Let $E\in (\gc)^s\cap i\g$ be an arbitrary non-trivial nilpotent element.
There exists an ${\mathfrak s\mathfrak l}_2$-triple $(E,H,F)$ in
$({\mathfrak g}^{\C})^{s},$ i.e., $[E,F]=H,\ [H,E]=2E$ and $[H,F]=-2F$
such that  $E,F\in {\mathfrak q}$ and $H\in {\mathfrak h}$.

\end {lem}

\begin {proof}

Since $(\gc)^s$ is reductive, there exists  a ${\mathfrak s\mathfrak l}_2$-triple
$(E,H,F)$  in
$({\mathfrak g}^{\C})^s$  by the theorem 
of  Jacobson-Morozov.  It can be chosen to be $\sigma$-compatible.
\s
To see this, split 
$H=H^\sigma+H^{-\sigma}$ with respect to the  $\sigma $-eigenspace
decomposition of $({\mathfrak g}^{\C})^s$. Since $[H,E]=2E$ and 
$\sigma(E)=-E$, it follows that $[H^{-\sigma},E]=0$. 
Hence, we may assume that $H=H^\sigma$
(see \cite {Bou} Chpt.~VIII, $\S$11, Lemme 6). 
Observe further that in this case
$[E,F]=[E,(F)^{-\sigma}]=H$ and $[H, (F)^{-\sigma}]=(F)^{-\sigma}$.
The desired result follows then from the uniqueness of the 
third element $F$ in a $\sL_2$-triple.
\end {proof} 
\s
Now we have all the ingredients which are needed 
to give a complete characterization of the closed orbits in $Im(\eta):$

\begin {prop} \label {closed orbits}

{\bf (Closed orbits)} If $\eta(x)=us$ is the Jordan decomposition,
then the orbit $G.\eta (x)=G.(su)$ 
contains the closed orbit $G.s$ in its closure $\overline{G.\eta (x)}$.
In particular, $G.\eta (x)$
is closed if and only if
$\eta (x)$ is semi-simple and
$s\in Im(\eta )$.

\end {prop}

\begin {proof}

Let $u=\Ad(\exp N)$ with $N$ as in Prop. \ref{lifting}. 
Hence, by Lemma \ref{sltriple} there is a ${\mathfrak s\mathfrak l}_2$-triple
$(N,H,F)$ ($E=N$) such that $[tH,N]=2tN$, i.e., $\Ad(\exp tH)(N)=e^{2t}N$ 
for every $t\in \R$. Note also that $\exp (\R H)\subset G\cap \exp (\gc)^s$
by construction of the ${\mathfrak s\mathfrak l}_2$-triple. It follows
that
\begin{eqnarray*}
 \eta(\exp tH\cdot x) & = & \exp tH.(us)=\Ad(\exp tH)\cdot us\cdot \Ad(\exp -tH)= \\
 & = & \Ad(\exp tH)\cdot u\cdot \Ad(\exp -tH) \cdot s= \\
 & = & \Ad(\exp tH)\Ad(\exp N) \Ad(\exp -tH) \cdot s= \\
  & = & \Ad(\exp e^{2t}N) \cdot s.
\end{eqnarray*}
For $t\to -\infty$ it follows $\exp tH.(us)=\Ad(\exp e^{2t}N)\cdot s\to s.$ 
Hence, the closed  orbit $G.s$ lies in the closure of $G.(us).$ 
In particular $G.(us)$ is non-closed if $u\ne 1$, i.e., if $\eta(x)$ is 
not semisimple. This, together with 
Lemma  \ref{ss then closed} implies that $G.\eta (x)$
is closed if and only if $\eta (x)$ is semi-simple.
Recall that the image $Im(\eta)$ is closed. This forces 
$s\in Im(\eta)$ and the proof is now complete.
\end {proof}

\subsection {Elliptic elements and closed orbits} \label {III.3}

Every non-zero complex number $z$ has the unique
decomposition $r\cdot e^{i\phi}$ into the hyperbolic part $r>0$ and
elliptic part $e^{i\phi}$.
This generalizes for an
arbitrary semisimple element $s\in \GL(\gc)$: 
By decomposing its eigenvalues one obtains the unique decomposition 
$s=s_{\rm ell}s_{\rm hyp}=s_{\rm hyp}s_{\rm ell}$.
An element $x\in G^{\C}$ is said to be {\it elliptic\/} if and
only if $\eta (x)=s$ is semi-simple with eigenvalues lying in
the unit circle.  It should be remarked that
$x$ itself may in such a case not be a semisimple element of the
group $G^\C$,e.g., $K^\C $ contains unipotent elements.
\s
Let $\Omega _{\rm ell}\subset G^\C $ be the set of elliptic elements.
This set is invariant by the right-action of $K^{\C}$, 
and therefore by choosing the same base point $x_0$ as in the case of 
$\OAG$, by abuse of notation we also regard $\Omega _{\rm ell}$ as a subset of
$G^{\C}/K^{\C}$. We reiterate that, since the map $\eta$ is 
not a group morphism, the classical notion of an elliptic element 
in $G^\C $ differs from the above definition.

\begin{lem} 

For $U$ the maximal compact subgroup of $G^{\C}$
defined by $\theta $ it follows that
$U\subset \Omega _{\rm ell}$.

\end{lem}

\begin{proof}

For $\theta $ the Cartan involution defining ${\mathfrak u}$, observe 
that $\widehat U:=\{ \varphi \in \Aut_{\R}({\mathfrak g}^{\C}):
\varphi \theta =\theta \varphi \} $ is a maximal compact subgroup of 
$\Aut_{\R}({\mathfrak g}^{\C})$ (with identity component 
$\Ad(U)$).
\s
Now $\theta $ commutes with every term in the definition of 
$\eta (u)$for every $u\in U$. It follows that 
$\theta \eta(u)=\theta \eta(u)$.  
Therefore $\eta (U)$ is contained in the compact group 
$\widehat U$ and consequently $U\subset \Omega _{\rm ell}$.  
\end{proof}
\s
\begin{prop}  \label{ell}
{\bf (Elliptic elements)} 
In the homogeneous space $G^{\C}/K^{\C}$ the
set of elliptic elements is described as
$\Omega _{ell}=G\cdot\exp(i{\mathfrak a})\cdot x_0$.

\end{prop}

\begin{proof}

Observe that $\Omega _{ell}$ is $G$-invariant. Hence, the above 
Lemma implies that
$G.\exp(i{\mathfrak a})\cdot x_0\subset \Omega _{ell}$.
\s
Conversely, suppose $x$ is elliptic, i.e., $\eta(x)$ is 
contained in some maximal compact subgroup of $\Aut_\R(\gc)$.
Hence, there is a Cartan involution $\theta '' :{\mathfrak g}^{\C}\to 
{\mathfrak g}^{\C}$ which commutes with $\eta (x)$. We
now make the usual adjustments so that, after replacing $x$ by
an appropriate  $G$-translate, $\eta (x)$ will commute with the given Cartan
involution $\theta $.
\s
For this, if $\theta ''$ does not commute with $\sigma $,  define
the semisimple element $\rho :=\sigma \theta ''\sigma \theta ''$ which
is diagonalizable with all positive eigenvalues over ${\R}$.  It follows that  
$\rho ^t$ is defined for all $t\in {\R}$, and 
$\theta ':=\rho ^{\frac 1 4}\theta ''\rho ^{-{\frac 1 4}}$ commutes 
with $\sigma $ (see \cite {Hel}, Chp. III, $\S 7$). By direct 
calculation one verifies that $\rho $, hence $\rho^t $ commutes 
with $\eta (x)$. Thus it follows that $\theta '$ and $\eta (x)$ commute.
\m
Finally, since $\theta '$ and our original $\theta $ both commute 
with $\sigma $, there exists $h\in G$ such that 
$\Ad(h)\theta '\Ad(h^{-1})=\theta $. Consequently, if $x$ is replaced 
by $h^{-1}\cdot x$, then we may assume that $\eta (x)$ and $\theta $ commute.
\m
Now we will adjust $x$ so that it lies in $U$.  With respect to the global 
Cartan decomposition  of $G^{\C}$
defined by $\theta $ 
write $x=u\exp(Z)$,  i.e., $u\in U$ and $\theta (Z)=-Z$. 
We now show that in fact $\exp(Z)\in K^{\C}$.
\s
Since $\theta $ commutes with
$\sigma $, $\tau $ and $u$ and anti-commutes with $Z$, we
have
\begin{eqnarray*}
\theta \eta (x)& = & 
\theta\cdot (\sigma \Ad(u)\Ad(\exp(Z))\tau \Ad(\exp(-Z))\Ad(u^{-1}))=  \\
 & = & \sigma \Ad(u)\Ad(\exp(-Z))\tau \Ad(\exp(Z))\Ad(u^{-1})\cdot \theta . 
\end{eqnarray*}
On the other hand
$$ \theta \eta(x)=\eta(x)\theta=\sigma \Ad(u)\Ad(\exp(Z))\tau 
\Ad(\exp(-Z))\Ad(u^{-1})\cdot \theta  .
$$
Combining these two equations, we obtain $\Ad(\exp(Z))\tau 
\Ad(\exp(-Z)=\Ad(\exp(-Z))\tau \Ad(\exp(Z))$
and consequently $\Ad(\exp(2Z)$ commutes with $\tau $.  Since
the restriction $\Ad:\exp(i{\mathfrak u})\to \Aut({\mathfrak g}^{\C})$ is 
injective, it follows that $\tau (\exp(Z))=\exp(Z)$, i.e., 
$\exp(Z)\in K^{\C}$. Replacing $x$ by $x\exp(-Z)$,
it follows that $x\cdot x_0=x\exp(-Z)\cdot x_0;$
hence, we may assume that $x\in U$.
\m
Since $U=K\cdot \exp(i{\mathfrak a})\cdot K$, we may 
assume that $x\in K\exp(i{\mathfrak a})$ and then translate it
by left multiplication by an element of $K$ to reach the following
conclusion: If $x\in G^\C$ is elliptic, then there exists $h\in G$ and 
$l\in K^\C$
with $hxl\in \exp i\a$ or, equivalently, there is $h\in G$ with
$hx\cdot x_0\in \exp(i{\mathfrak a})\cdot x_0$.  This proves the
inclusion $\Omega _{\rm ell}\subset G\cdot \exp(i{\mathfrak a})\cdot x_0$.
\end {proof}
\s
The above proposition yields key information on the set 
$\Omega _{c\ell }$ of closed $G$-orbits in
$G^{\C}/K^{\C}$. For a subset $M$ of $\Omega $
let $\text{\rm c}\ell (M)$ denote its 
topological closure.

\begin{cor} \label{elliptic}

$$\Omega _{c\ell}\cap \text{\rm c}\ell(\OAG)=
G\cdot \text{\rm c}\ell(\exp(i\oag )x_0)=
\Omega _{\rm ell}\cap \text{\rm c}\ell(\OAG).$$

\end{cor}

\begin{proof}

By the above Proposition, $\OAG=G\cdot \exp(i\oag )\cdot x_0
\subset \Omega _{ell}$. Thus, by continuity, if $x\in \text{\rm c}\ell(\OAG)$, 
then the semisimple part of $\eta(x)$ is elliptic.  As a result, if
$x\in \Omega _{c\ell}\cap \text{\rm c}\ell(\OAG)$, then $\eta(x)$
is semisimple by Prop \ref{closed orbits} and hence it is elliptic. It follows
that $x\in G\cdot \exp(i{\mathfrak a})\cdot x_0$. 
But $\exp(i{\mathfrak a})\cap \text{\rm c}\ell(\OAG)=
\text{\rm c}\ell(\exp(i\oag ))$ 
(``$\supset$'' is obvious and ``$\subset$'' follows
by transversality of $G$--orbits at $\exp (i\a)\cdot x_0$) and therefore
$\Omega _{c\ell}\cap \text{\rm c}\ell(\OAG)\subset
G\cdot \text{\rm c}\ell(\exp(i\oag ))$.
\s
As a consequence of the above Proposition and Lemma \ref{ss then closed} 
we see that
$G\cdot \text{\rm c}\ell(\exp(i\oag ))\subset
\Omega _{\rm ell}\cap \text{\rm c}\ell(\OAG)$.
\s
Finally, if $x\in \Omega _{ell}$, then in particular
$\eta (x)$ is semisimple, and therefore $G\cdot x$ is closed.
Hence,
$\Omega _{ell}\cap \text{\rm c}\ell(\OAG)
\subset \Omega _{c\ell}\cap \text{\rm c}\ell(\OAG)$.
\end{proof}
\s

\section {$Q_2$-slices}  \label {IV}

At a generic  point $y\in {\rm bd}(\OAG )$ we determine a
$3$-dimensional, $\sigma $-invariant, semi-simple subgroup $S^{\C}$
such that $S=(S^\C)^\sigma={\rm Fix}(\sigma :S^\C\to S^\C)$ is a non-compact real form 
and such that the isotropy group $S^{\C}_y$ is either a maximal 
complex torus or its normalizer. Geometrically speaking, 
$Q_2=S^{\C}\cdot y$ is either the $2$-dimensional affine quadric, 
which can be realized by the diagonal action as the complement of the diagonal 
in $\P_1({\C})\times\P_1({\C})$, or
its (2-1)-quotient, which is defined by exchanging the factors and
which can be realized as the complement of the (closed) $1$-dimensional
orbit of $SO_{3}(\C)$ in $\P_2({\C})$. By abuse of notation, we refer in both cases to
$S^{\C}\cdot y$ as a 2-dimensional affine quadric.
\s
The key property is
that, up to the above mentioned possibility of a (2-1)-cover,
the intersection $Q_2\cap \OAG$ is the Akhiezer-Gindikin
domain in $S^{\C}/K_S^{\C}$ for the unit disk
$S/K_S$.
\s
For the sake of brevity we say that the orbit $S^\C\cdot y$
is a $Q_2$-slice at $y$ whenever it has all of the above properties.

\subsection {Existence}  \label {IV.1}

Given a non-closed $G$-orbit $G\cdot y$ in $\text{\rm bd}(\OAG)$, we may
apply Prop. \ref {lifting} to obtain a lifting of the semi-simple
(elliptic) part of the Jordan decomposition of $\eta (x) $. For an
appropriate base point $z$ this lifting can be
chosen in $\text{\rm bd}(\exp(i\oag )$. Recall that the action
$G^\C\times \Aut_\R(\gc)\to  \Aut_\R(\gc)$ is given by conjugation
(see 3.1).
Note that the isotropy Lie algebra at $\varphi\in \Aut_\R(\gc)$ is 
the totally real
subalgebra of fixed points $(\gc)^\varphi=\{Z\in \gc:\varphi(Z)=Z\}.$

\begin {lem}  \label {optimal}

{\bf (Optimal base point)} Every non-closed $G$-orbit $G\cdot y$ in 
$\text{\rm bd}(\OAG)$ contains a point $z=\exp E\cdot \exp iA \cdot x_0$
such that $E\in (\gc)^{\eta(\exp iA)}\cap i\g $ is a non-trivial
nilpotent element.

\end{lem}

\begin{proof}

Let $\eta(y)=su$ be the Jordan decomposition and let
$N\in (\gc)^s\cap i\g$ be as in Prop. \ref {lifting}. We then have 
$$\textstyle \eta(\exp y)=
\eta(\exp(-\frac{1}{2}N)\exp(\frac{1}{2}N)\cdot y)=
\Ad(\exp N)\circ\eta(\exp(\frac{1}{2}N)\cdot y)=u\cdot s.$$
By Prop. \ref{closed orbits} and Cor. \ref{elliptic} the 
semisimple element $ \eta(\exp(\frac{1}{2}N)\cdot y)$ is elliptic. 
Hence, Prop. \ref{ell}
implies the existence of  $g\in G$ and $A\in {\rm bd}(\oag)$ 
such that $\exp\frac{1}{2}N\cdot y=g^{-1}\exp iA\cdot x_0$.
\s
Define now $E:=\Ad(g)(-\frac{1}{2}N)$ and observe that $g\cdot y=\exp E\exp iA\cdot x_0.$
Finally, $E\in (\gc)^{g.s}=(\gc)^{\eta(\exp iA)}$, and the lemma is proved.
\end{proof}
\m
Recall that $(\gc)^{\eta(\exp iA)}$ is a $\sigma$-stable real 
reductive algebra. Let $(\gc)^{\eta(\exp iA)}=\h\oplus \q$ be the 
decomposition into $\sigma$-eigenspaces. 
In this notation, the nilpotent element $E$ as in the above 
lemma belongs to $\q$.
\m
Let now an arbitrary  non-closed orbit $G\exp E\exp iA\cdot x_0$ be given.
Fix  a  ${\mathfrak s\mathfrak l}_2$-triple $(E,H,F)$
as in Lemma  \ref {sltriple}.
Let $S^{\C}$ be the complex subgroup of $G^{\C}$
defined by this triple. Set
$e:=iE$, $f:=-iF$ and let $S$ be the $\sigma $-invariant
real form in $S^\C.$ The Lie algebra of $S$ is then the subalgebra 
generated by the $\sL_2$ triple $(e,H,f)$.  Finally,
let $x_1=\exp(iA)\cdot x_0$ be the base point chosen as above
in the closure of a given $G$-orbit.

\begin {lem}  \label {Sisotropy}

The connected component $(S^{\C}_{x_1})^\circ $ of the
$S^{\C}$-isotropy at $x_1$ is the $1$-parameter subgroup
$\{ \exp(zH):z\in \C \}\cong \C^* $.

\end {lem}

\begin {proof}

Since the action $S^\C\times \GCKC\to \GCKC$ is affine-algebraic, the
orbit $S^\C\cdot \exp iA\cdot x_0=S^\C\cdot x_1$ is an affine variety. 
Then the isotropy
at $\exp iA\cdot x_0$ is 1-dimensional or $S^\C.$ Note that 
$S^\C\cdot x_1$ cannot be a point, because by construction 
$\exp E \cdot x_1\ne x_1$; therefore $S^{\C}_{x_1}$
is 1-dimensional.
\s
We now show that $\exp \R H\cdot x_1=x_1$, or equivalently, 
$\exp tH.\eta(x_1)=\eta(x_1)$. 
Define $\varphi :=\Ad(\exp(iA))\tau \Ad(\exp(-iA))$ and note that
$H \in (\gc)^\sigma\cap (\gc)^\varphi=\h$ yields
\begin{eqnarray*}
\exp tH.\eta(x_1)&=&\Ad(\exp tH)
\circ \sigma\varphi\circ \Ad(\exp -tH)= \\
&=& \Ad(\exp tH)\Ad(\exp -tH)\circ \sigma\varphi=\eta(x_1).
\end{eqnarray*}
It follows that $\exp \C H\cdot x_1=x_1$. Since $S^{\C}_{x_1}$ 
is 1-dimensional and $H$ semisimple,
we deduce $(S^{\C}_{x_1})^\circ=\exp \C H\cong \C^*$.
\end {proof}

\subsection {Genericity} \label {IV.2}

Without going into a technical analysis of $\text{\rm bd}(\OAG )$, we will
construct $Q_2$-slices only at its generic points.
The purpose of this section is to introduce the appropriate
notion of "generic" and prove that the set of such points
is open and dense. The set of generic points is defined to be the
complement of the union of small semi-algebraic sets
${\cal C}$ and ${\cal E}$ in $\text{\rm bd}(\OAG )$.  We begin
with the definition of ${\cal C}$.
\m
Let $R:=\text{\rm bd}(\exp(i\oag )\cdot x_0)$ and recall that for
$y\in \text{\rm bd}(\OAG )$ the orbit $G\cdot y$ is closed if and only
if $G\cdot y\cap R\ne\varnothing $.  In fact $R$ 
parameterizes the closed orbits in $\text{\rm bd}(\OAG )$ 
up to the orbits of afinite group.  Recall
also that $R$ is naturally identified with $\text{\rm bd}(\oag )$,
which is the boundary of a convex polytope, and is defined by 
linear inequalities.  Let $E$ be the
image in $R$ of the lower-dimensional edges in $\text{\rm bd}(\oag )$,
i.e., the set of points which are contained in at least two
root hyperplanes $\{\alpha=c_\alpha\}$.  
Finally, let $R_{\text{\rm gen}}:=R\smallsetminus E$.  
\m
As we have seen in Cor. \ref{elliptic}, the set of closed 
orbits in the boundary of $\OAG$ can be described as  
$G\cdot\text{\rm bd}(\exp(i\oag)x_0)$.
This is by definition the set ${\cal C}$.

\begin{lem}

For $x\in R$ it follows that $dim\, G.x\le \codim_\Omega {\rm bd}(\Omega) -2$.

\end{lem}

\begin {proof}
Note that ${\rm bd}(\Omega)$ is connected and of codimension 1 in $\Omega$.
The $G$-isotropy group $C_K({\mathfrak a})$ at generic points of 
$exp(i\omega _{AG}).x_0$ fixes this slice pointwise and therefore
is contained in a maximal compact subgroup of the isotropy
subgroup $G_x$ of each of its boundary points. Since by definition
$G_x$ is non-compact, it follows that $dim\ G_x$ is larger than
the dimension of the generic $G$-isotropy subgroup at points
of $exp(i\omega _{AG}).x_0$. 
\end{proof}
\s
\begin {rem}

For $x\in {\rm bd}(\Omega)_{{\rm gen}}$ the isotropy 
subgroup $G_x$ is precisely calculated in $\S$\ref {IV.3}. This shows that 
$\dim G\cdot x=\codim_\Omega {\rm bd}(\Omega)-m$, where $m$ is at least 2.
Thus, by semi-continuity
we have the estimate $\dim\, G.x\le \codim_\Omega {\rm bd}(\Omega)-m$ for all
$x\in {\rm bd}(\Omega)$.
\end {rem}
\s
Now let $X:=Im(\eta )\subset Aut_\R({\mathfrak g}^\C)$.
It is a connected component of a real  algebraic submanifold in $\Aut_\R(\gc)$.
The complexification $X^\C$ of $X$ which is contained in the complexification
$\Aut_\C({\mathfrak g}^\C\times {\mathfrak g}^\C)$
of $\Aut_R({\mathfrak g}^\C)$ is biholomorphic to 
$G^\C/N^\C\times G^\C/N^\C$, where $N^\C$ denotes the normalizer of 
$K^\C$ in $G^\C$.  The complexification
of the piecewise real analytic variety $R$ is a piecewise
complex analytic subvariety $R^\C$ of $X^\C$ defined 
in a neighborhood of $R$ in $X^\C$.  Finally, 
let $\pi :X^\C\to X^\C/\!/G^\C$
be the complex categorical quotient.
\s
Recall that in every $\pi $-fiber there is a unique closed 
$G^\C$-orbit.  The closed $G$-orbits in $X$ are components
of the the real points of the closed $G^\C$-orbits which
are defined over $\R$. For a more extensive discussion of the interplay betweenthe real and complex points in complex varieties defined over $\R$ see 
\cite {Sch1}, \cite {Sch2}, \cite {Br}.
\m
Let $k:=\dim_\R\Omega -\dim\, R-m$ be the dimension of the
generic $G$-orbits of points of $R$ and let 
$S_k$ be the closure in $X^\C$ of
$\{ z\in X^\C:G^\C\cdot z \ \text{\rm is closed and} 
\ k\text{\rm -dimensional}\} $. Define
${\cal C}_k:=S_k\cap R^\C$. It follows that
${\cal C}_k$ is a piecewise complex analytic set of dimension
$k+\dim_\C R^\C$.

\begin {prop}

The set 
$G\cdot R=\{ x\in \text{\rm bd}(\OAG ): G\cdot x \ \text{\rm is closed}\} $
is contained in a closed semi-algebraic subset ${\cal C}$ of codimension 
at least one in $\text {\rm bd}(\OAG )$.

\end {prop}

\begin {proof}

The set ${\cal C}$ is defined to be the intersection of the real
points of ${\cal C}_k$ with $\text{\rm bd}(\OAG )$.  The desired
result follows from
$\dim _\C{\cal C}_k = k+\dim_\C R^\C$.
\end {proof}
\m
Recall that $\pi$ denotes the categorical quotient map $\pi:X^\C\to X^\C/\!/ G^\C$.
Define ${\cal E}:=\eta ^{-1}(\pi ^{-1}(\pi (E)))\cap 
\text{\rm bd}(\OAG )$. In particular it is a closed
semi-algebraic subset of $\text{\rm bd}(\OAG )$ which
contains the set $\{x\in \text{\rm bd}(\OAG ):
\text{\rm c}\ell(G\cdot x)\cap E\not=\emptyset \} $. 

\begin {defn}

A point $z\in \text{\rm bd}(\OAG )$ is said to be generic if it is
contained in the complement of ${\cal C}\cup {\cal E}$.
\end {defn}
\s
Let $\text{\rm bd}_{\text{\rm gen}}(\OAG )$ denote the set of generic
boundary points.

\begin {prop} \label {generic}

The set of generic points $\text{\rm bd}_{\text{\rm gen}}(\OAG )$ is
open and dense in $\text{\rm bd}(\OAG )$.
\end {prop}
\s
It has already been noted that ${\cal C}$ and ${\cal E}$ are closed.
Since ${\cal C}$ is of codimension two, the complement of ${\cal C}$ is dense.
Thus this proposition is an immediate consequence of the following fact.

\begin {prop} \label {codimension two}

The saturation ${\cal E}$ is at least $1$-codimensional
in $\text{\rm bd}(\OAG )$.
\end {prop}
\s
This in turn follows from a computation of the dimension of the
fibers at points of $E$ of the above mentioned categorical quotient. 
For this it is convenient to use the Jordan decomposition 
$\eta (z)=u\cdot s$ for $z\in \Omega $ such that
$x=\exp(i{\mathfrak a})\cdot x_0$ is in $\text{\rm c}\ell(G\cdot z)$.
\m
As in Lemma \ref {optimal} we choose an optimal base point such that
$\eta (x)=s$ and $u=\Ad(\exp(N))$ with $N\in {\mathfrak q}$,
where ${\mathfrak h}\oplus {\mathfrak q}$ is the $\sigma $-decomposition
of ${\mathfrak l}=({\mathfrak g}^\C)^s$.  Let $N_x$ be
the the cone of nilpotent elements in ${\mathfrak q}$ and observe
that the saturation ${\cal E}_x=\{z\in \text{\rm bd}(\OAG ):
x\in \text{\rm c}\ell(G\cdot z)\} $ is an $N_x$-bundle over the closed
orbit $G\cdot x$. Thus it is necessary to estimate
$dim_\R N_x$.   
\m
Recall that any two maximal toral Abelian subalgebras of 
${\mathfrak q}^\C$ are conjugate and therefore the dimension $m$
of one such is an invariant.  Since ${\mathfrak a}^\C$
is such an algebra, the following is quite useful 
(see \cite {KR}).

\begin {lem} 

The complex codimension in ${\mathfrak q}^\C$
of every component of the nilpotent cone in 
${\mathfrak q}^\C$ is $m$.

\end {lem}
\s
{\it Proof of Prop. \ref {codimension two}}
\s
We prove the estimate $\codim_\Omega {\cal E}_x\ge \dim\ {\mathfrak a}$.
For this observe that, since $G.s$ is closed in 
$Aut_\R({\mathfrak g}^\C)$, an application of the
Luna slice theorem for the (closed) complex orbit $G^\C$.s
in the complexification of $Im(\eta )$ yields the bundle structure 
$Im(\eta )=G\times _{G_s}{\mathfrak q}$ locally near $s$; in particular
$\codim _{\mathfrak q}N_x=\codim _\Omega ({\cal E}_x)$. The result follows
from the above Lemma by noting that 
$\codim_{\mathfrak q}N_x$ is at most the complex codimension of the nilpotent
cone in ${\mathfrak q}^\C$ and, as mentioned above,
that $\dim  {\mathfrak a}=m$.  $\square $  
\m
The group $S^{\C}$ constructed above for a generic boundary point 
has the property that the intersection
of the $S^{\C}$-orbit, i.e., a 2-dimensional affine 
quadric $Q_2\cong SL_2(\C)/\C^*$ (or $\cong SL_2(\C)/N(\C^*)$) 
with  $\OAG$ contains an Akhiezer-Gindikin domain 
$\Omega_{\rm AG}^{SL}\cong D\times \overline D$ of $Q_2$. 
To see this, we will conjugate $S^\C$ by an element
of $G$ in order to relate $S^\C$ to the fixed Abelian Lie algebra $\a$.
This is carried out in the next section.

\subsection {The intersection property} \label {IV.3}

To complete our task we conjugate the group $S^{\C}$
obtained in $\S$\ref {IV.1} above by an element 
$h$ in the isotropy group  $G_{x_1}$ so that it can be 
easily seen that the resulting
orbit $Q_2=S^{\C}\cdot x_1$ intersects $\OAG$
in the Akhiezer-Gindikin domain of $Q_2$.
\s
The following is a first step in this direction.

\begin {prop}  \label {conjugate}

Let $G\exp E\exp iA\cdot x_0=G\cdot x_1$ be any non-closed 
orbit in ${\rm bd}(\OAG)$ and 
$(E,H,F)$ a ${\mathfrak s\mathfrak l}_2$-triple in 
$(\gc)^{\eta(\exp iA)}$ as in Lemma \ref{sltriple}. 
Given  $Z:=E-F$, then there exits
$h\in G_{x_1}$ so that $\Ad(h)(Z)\in i{\mathfrak a}$.
\end {prop}
\s
This result is an immediate consequence of the following basic
fact.

\begin {lem} \label {conjinq}

Let ${\mathfrak l}$ be a real reductive Lie algebra, $\theta $
a Cartan involution and $\sigma $ a further involution which
commutes with $\theta $.  Let 
${\mathfrak l}={\mathfrak k}\oplus {\mathfrak p}$ be the 
eigenspace decomposition with respect to $\theta $ and
${\mathfrak l}={\mathfrak h}\oplus {\mathfrak q}$ with respect
to $\sigma $.  Then, if ${\mathfrak a}\subset {\mathfrak p}\cap {\mathfrak q}$
is a maximal Abelian subalgebra of ${\mathfrak q}$ and
$\xi $ is a hyperbolic semisimple element of ${\mathfrak q}$, there
exists $h\in Int({\mathfrak h})$ such that 
$\Ad(h)(\xi )\in {\mathfrak q}$.

\end {lem}

\begin {proof}

Since $\xi $ is hyperbolic, we may assume that there is a Cartan involution
$\theta ':{\mathfrak l}\to {\mathfrak l}$ such that 
$\theta '(\xi )=-\xi $ and $\theta'\sigma=\sigma\theta'$.
Then there exists 
$h\in Int({\mathfrak h})$ with $\Ad(h)\theta '\Ad(h^{-1})=\theta $
(see \cite {M1}) and $\Ad(h)(\xi )\in {\mathfrak p}\cap {\mathfrak q}$.
\s
To complete the proof, just note that 
$({\mathfrak h}\cap {\mathfrak k})\oplus ({\mathfrak p}\oplus {\mathfrak q})$
is a Riemannian symmetric Lie algebra where any two maximal Abelian
algebras in ${\mathfrak p}\cap {\mathfrak q}$ are conjugate
by an element of $Int({\mathfrak h}\cap {\mathfrak k})$.
\end {proof}
\m
{\it Proof of Prop. \ref {conjugate}.}
Observe that ${\rm ad}(Z)$ has only imaginary eigenvalues.
Replacing $(\gc)^{\eta(\exp iA)}=\h\oplus\q $ by the dual 
$\tilde {\mathfrak l}:={\mathfrak h}\oplus i{\mathfrak q}=
\tilde {\mathfrak h}\oplus \tilde {\mathfrak q}$ and defining
$\tilde \sigma $ and $\tilde \theta $ accordingly, we apply the
above Lemma to $\xi :=iZ$ and the Abelian Lie algebra 
${\mathfrak a}\subset \tilde{\mathfrak q}$ to obtain 
$h\in Int({\mathfrak h})$ with $\Ad(h)(\xi )\in {\mathfrak a}$.
Thus $\Ad(h)(Z)$ has the required property $\Ad(h)(Z) \in i {\mathfrak a}$. \hfill $\square $ 
\b
We now show that for $z\in \text{\rm bd}_{\text {\rm gen}}(\OAG )$
the group $S^\C$ which is associated to the 
${\mathfrak s\mathfrak l}_2$-triple constructed in the above
proposition as produces a $Q_2$-slice.  For a precise formulation 
it is convenient to let $\text{\rm bd}_{\text{\rm gen}}(\oag ):=
\text{\rm bd}(\oag )\smallsetminus E$, where $E$ is the union of
the lower-dimensional strata as in $\S$\ref {IV.2}.

\begin {prop}

For $z\in \text{\rm bd}_{\text{\rm gen}}(\OAG )$ and 
$x_1=\exp iA\cdot x_0$ the associated point
with $iA\in \text{\rm bd}_{\text{\rm gen}}(i\oag )$
it follows that the 
line ${\R}(E-F)$ is transversal to $\text {\rm bd}_{\text{\rm gen}}(i\oag)$ 
at $iA$ in $i\a$.
\end {prop}
\s
The proof requires a more explicit description of 
$({\mathfrak g}^\C)^{\eta (\exp iA)}=:\l=\h\oplus \q$.
For this, recall the root decompositions of $\g$ and $\gc$  
with respect to $\a$ or $\ac,$ respectively:
$ \gc=C_{\kc}(\ac)\oplus \ac\oplus\bigoplus_{\Phi(\a)}\gc_\lambda$. 
The behavior of this decomposition with respect to our involutions 
is the following: 
$\theta(\gc_\lambda)=\gc_{-\lambda}$ and  
$\tau(\gc_\lambda)=\gc_{-\lambda}$;
furthermore, the root decomposition is $\sigma$--stable, 
i.e., $\sigma(\gc_\lambda)=\gc_{\lambda}$.
Fix a $\tau$--stable basis of root covectors, i.e., select 
any basis $L_\lambda^1,...\,,L_\lambda^k$ of 
$\g_\lambda=(\gc_\lambda)^\sigma$and define 
$L_{-\lambda}^j:=\tau(L_{\lambda}^j).$
Define $\g{[\lambda]}:=\g_\lambda\oplus 
\g_{-\lambda},\ \gc{[\lambda]}:=\g{[\lambda]}\oplus i\g{[\lambda]}$ 
and notice that $\g{[\lambda]}=(\g{[\lambda]})^\tau
\oplus (\g{[\lambda]})^{-\tau}.$ Finally, set
$$X_{[\lambda]}^j:=L_{\lambda}^j+L_{-\lambda}^j\qquad
\qquad Y_{[\lambda]}^j:=L_{\lambda}^j-L_{-\lambda}^j $$
and observe that $X_{[\lambda]}^j\in 
(\g{[\lambda]})^\tau, \ Y_{[\lambda]}^j\in (\g{[\lambda]})^{-\tau}.$ 
The reason for introducing such a basis is that
the complex subspaces $(\!(X_{[\lambda]}^j,Y_{[\lambda]}^j)\!)_\C$ are 
$\Ad(t)$-stable for any $t:=\exp iA,$ $A\in \a.$
\s
Express $\Ad(t)$ as a matrix  with respect to the basis 
$X_{[\lambda]}^j,Y_{[\lambda]}^j$:
$$\Ad(t)\big\vert_{\textstyle(\!(X_{[\lambda]},Y_{[\lambda]})\!)}=
\begin{pmatrix}
\cosh \lambda(iA) &\sinh \lambda(iA)\cr\sinh \lambda (iA)&\cosh \lambda(iA)\cr
\end{pmatrix}
$$
Let $\gc=\kc\oplus \pc$ be the complexification of the 
Cartan decomposition of $\g$.
A simple calculation yields for $t=\exp iA$:
\begin{eqnarray*}
\h\kern2pt=\kern 4pt\g\cap \Ad(t)(\kc)&=&
\displaystyle C_\k(\a)\oplus\bigoplus_{\lambda(A)=\Z\pi}\!\g[\lambda]^
\tau\oplus\kern-1em\bigoplus_{\lambda(A)={\frac \pi 2}+\Z\pi}\!\!\!\g
[\lambda]^{-\tau} \\
\q= i\g\cap \Ad(t)(\pc)&=&
\displaystyle i\a\oplus\bigoplus_{\lambda(A)=\Z\pi}\!i\g[\lambda]^{-\tau}\oplus
\bigoplus_{\lambda(A)={\frac \pi 2}+\Z\pi}\!\!\!i\g[\lambda]^\tau . 
\end{eqnarray*}
Let $A\in \text{\rm bd}_{\text{\rm gen}}(\oag )$ be boundary-generic,
i.e., there is a single $\lambda\in \Phi(\a)$ with
$$ \lambda(A)=\pm \frac{\pi}{ 2} \qquad\qquad 
\mu(A)\not \in \frac{\pi}{2}\Z \qquad\hbox{for all } 
\ \mu\in \Phi(\a)\smallsetminus \{\pm\lambda\}.$$
The above general formulas imply that  the
 centralizer subalgebra 
$(\gc)^{\eta(\exp iA)}$ for such a boundary-generic point as above 
is  given by
$$
(\gc)^{\eta(\exp iA)}=\mm\oplus \g[\lambda]^{-\tau} 
\oplus i\a \oplus i\g[\lambda]^{\tau}.
$$
To complete the proof of the proposition it is then enough to show that
for the selected $\sL_2$-triple $(E,H,F)\in \gc)^{\eta(\exp iA)}$ it 
follows that $E-F\in \R ih_\lambda$,
where $h_\lambda\in \a$ is the coroot determined by the root 
$\lambda\in \Phi(\a)$. This is the content of the following

\begin{lem} \label{h_lambda}

Let $A\in \{\lambda=\pi/2\}\cap {\rm bd}_{\rm gen}(\oag)$
be boundary-generic as above. Then $E-F\in \R ih_\lambda$.

\end{lem}

\begin{proof}

Let $\l:= (\gc)^{\eta(\exp iA)}=\h\oplus\q.$ Since $(\!(E,H,F)\!)_\R$ is 
semisimple, it follows that $(\!(E,H,F)\!)_\R\subset[\l\colon\l]$. Hence,
since $B([\g_\lambda:\g_{-\lambda}],\{\lambda=0\})=0$ 
($B$ denotes the Killing form) we have
\begin{eqnarray*}
[\l:\l]&=& [\mm\oplus \g[\lambda]^{-\tau} 
\oplus i\a \oplus i\g[\lambda]^{\tau}:\mm\oplus \g[\lambda]^{-\tau} 
\oplus i\a \oplus i\g[\lambda]^{\tau}]= \\
&=& \mm\oplus \R ih_\lambda\oplus \g[\lambda]^{-\tau}
\oplus i\g[\lambda]^{\tau}
\end{eqnarray*}
By Prop. \ref{conjugate} we have $E-F\in i\a.$ Finally, 
$E-F\in i\a\cap [\l:\l]=\R ih_\lambda$.
\end{proof}
\m
Recall that the set ${\rm bd}_{\rm gen}(\OAG )=
{\rm bd}(\OAG)\smallsetminus ({\cal C}\cup{\cal E})$ consists of 
certain non-closed orbits in the boundary of $\OAG$.

\begin {thm}  \label {slice}

On every $G$-orbit in $\text{\rm bd}_{\text{\rm gen}}(\OAG )$ there exists 
a point  of the form $z:=\exp E\exp iA\cdot x_0$, 
$A\in {\rm bd}_{\rm gen}(\oag )$, $E$ nilpotent,
 and a corresponding 3-dimensional simple subgroup $S^\C\subset G^\C$ such that
\begin{enumerate}

\item The 2-dimensional affine quadric $S^\C\cdot\exp iA\cdot x_0 =: 
S^\C\cdot x_1$  contains $z$
\item The intersection  $\OAG\cap S^\C\cdot x_1 $  contains an 
Akhiezer-Gindikin domain $\Omega_{\rm AG}^{SL}$ of
$S^\C\cdot x_1$, i.e., the orbit $S^\C\cdot x_1$ is a $Q_2$-slice.
\end{enumerate}

\end {thm}

\begin{proof}
Given a  non-closed $G$-orbit in $\text{\rm bd}_{\text{\rm gen}}(\OAG )$ 
let $z=\exp E\exp iA\cdot x_0$ an optimal base point as in 
Lemma \ref{optimal}. By Prop. \ref{conjugate} we may choose an 
$\sL_2$-triple $(E,H,F)$ in $(\gc)^{\eta(\exp iA)}$ such that $E-F\in i\a$.
Let $S^\C\subset G^\C$ be the complex subgroup with Lie algebra 
$\sc:=(\!(E,H,F)\!)_\C$. By construction $S^\C\cdot x_1\ni z$.
\s
For a boundary-generic point $x_1$ with $\lambda(A)=\pi/2$ 
and $\mu(A)\ne \frac{\pi}{2}\mathbb Z$ for all $\mu\ne \pm \lambda$ 
we already know by \ref{h_lambda} that $E-F\in \R ih_\lambda.$ 
Assume that $h_\lambda\in \a$ is the normalized coroot of $\lambda$, i.e., 
$\lambda(h_\lambda)=2$. Since $\oag$ is invariant under the 
Weyl group, the image $A^\prime$
of $A$ under the reflection on $\{\lambda=0\}$ is also boundary-generic, 
and the intersection of $A-\R h_\lambda$ with $\oag$ is the 
segment $\{A-th_\lambda:t\in (0,\pi/2)\}$ with boundary points 
$A$ and $A^\prime$. 
\s
Recall that $(e,H,f)$ with $E=ie$ and $F=-if$ 
is a $\sL_2$-triple in $\sc$ such that $\ss:=
\g\cap \sc=(\!(e,H,f)\!)_\R$. Let $S$ denote the
corresponding subgroup in $S^C$ (isomorphic to $SL_2(\R)$ or $PSL_2(\R)$).
The $S$-isotropy at all points $\exp ((-\frac{\pi}{2},0)ih_\lambda+iA)
\cdot x_0$ is compact and it is non-compact at 
$\exp iA\cdot x_0$ and $\exp iA'\cdot x_0$.
Hence, $S\cdot \exp((-\frac{\pi}{2},0)ih_\lambda+iA)\cdot x_0  $ is an
Akhiezer-Gindikin domain in $S^\C\cdot x_1$ which is contained in $\OAG$.
\end{proof}
\s

\subsection {Domains of holomorphy} \label {IV.4}

Let $S^{\C}=SL_2({\C})$, $S=SL_2({\R})$ be
embedded in $S^{\C}$ as the subgroup of matrices which
have real entries and let $K_S=SO_2(\R)$.  To fix the notation, 
let $D_0$ and $D_\infty$ be the open $S$-orbits in $\P_1({\C})$. 
Further, choose $\C\subset \C\P^1=\C\cup\{\infty\}$
in such a way that  $0\in D_0$ and
$\infty \in D_\infty$ are the $K_S$-fixed points.
\m
Now let $S^{\C}$ act diagonally on 
$Z=\C\P^1\times \C\P^1$ and
note that the open orbit $\Omega $, which is the complement of the
diagonal $\diag (\C\P^1)$ in $Z$, is the complex symmetric space 
$S^{\C}/K^{\C}_S$. 
Note that in $\C\P^1\times \C\P^1$ there are 4 open 
$SL_2(\R)\times SL_2(\R)$--orbits: the bi-disks $D_\alpha \times 
D_\beta$ for any pair $(\alpha,\beta)$ from $\{0,\infty\}$. 
As $S$-spaces, the domains $D_0\times D_\infty$ and $D_\infty\times D_0$
are equivariantly biholomorphic; further, they are actually subsets 
of $\Omega $,and the Riemannian symmetric space $S/K_S$ sits in
each of them as the totally real $S$-orbit 
$S\cdot(0,\infty)$ (or $S\cdot (\infty,0)$, respectively).
Depending on which of these points is chosen as a reference point 
in $\Omega$, both domains can be considered as the Akhiezer-Gindikin domain
$$
\OAG=D_0\times D_\infty=S\cdot \exp i\oag\cdot (0,\infty)\qquad
D_\infty\times D_0=S\cdot \exp i\oag\cdot (\infty,0)
$$
with $\oag=(-\frac{\pi}{4},\frac{\pi}{4})h_\alpha$ and 
$h_\alpha\in \a$ is the normalized coroot (i.e.,$\alpha(h_\alpha)=2$).
\m
Our main point here is to understand $S$-invariant Stein domains
in $\Omega $ which properly contain $\OAG $. By symmetry we may assume
that such has non-empty intersection with $D_0\times D_0$. 
Observe that 
$(D_0\times D_0) \cap \Omega=D_0\times D_0\smallsetminus \diag(D_0) .$
Furthermore, other than $\diag(D_0)$, all $S$-orbits in 
$D_0\times D_0$ are closed real hypersurfaces. For
$D_0\times D_0\smallsetminus \diag(D_0)$ let $\Omega (p)$ be the domain
bounded by $S\cdot p$ and $\diag(D_0)$. We shall show that a function
which is holomorphic in a neighborhood of $S\cdot p$ extends holomorphically
to $\Omega (p)$.
\m
For this, define $\Sigma :=\{ (-s,s):0\le s<1\} \subset D_0\times D_0$.
It is a geometric slice for the $S$-action.
We say that a ($1$-dimensional) complex curve 
$C\subset {\C}^2 \subset Z$ is a supporting curve for
${\rm bd }(\Omega(p))$ at $p$ if $C\cap \text{\rm c}\ell(\Omega(p))=\{ p\} $.
Here, $\text{\rm c}\ell(\Omega(p))$ denotes the topological closure in 
$D_0\times D_0$.

\begin{prop}

For every $p \in D_0\times D_0\smallsetminus \diag (D_0)$ there exists
a supporting curve for ${\rm bd}(\Omega(p))$ at $p$.

\end{prop}

\begin {proof}

Recall that we consider $D_0$ embedded in $\C$ as the unit disc.
It is enough to construct such  a curve $C\subset \C^2$ at each point
$p_s=(-s,s)\in \Sigma $, $s\ne 0$.  For this we define
$C_s:=\{ (-s+z,s+z):z\in {\C}\} $.  To prove 
$C_s\cap \text{\rm c}\ell(\Omega(p_s))=\{ p_s\}$
let $d$ be the Poincare metric of the unit disc $D_0$, considered 
as the function $d:D_0\times D_0\to \R_{\ge 0}$.
Note that it is an $S$-invariant function on $D_0\times D_0$.
In fact the values of $d$ parameterize the $S$-orbits.
\m
We now claim that $d(-s+z,s+z)\ge d(-s,s)=d(p_s)$ for $z\in \C$ 
and $(-s+z,s+z)\in D_0\times D_0$,
with equality only for $z=0$, i.e., $C_s$ touches 
$\text{\rm c}\ell( \Omega(p_s))$ only at $p_s$.
To prove the above inequality, it is convenient to compare
the Poincare length of the Euclidean
segment ${\rm seg}(z-s,z+s)$ in $D_0$ with the length of ${\rm seg}(-s,s)$.
Writing the corresponding integral for the length, it is clear,
 without explicit calculation, that $d(-s+x,s+x)>d(-s,s)$ for 
$z=x\in {\R}\smallsetminus 0$. The same argument shows also that
$d(-s+x+iy,s+x+iy)> d(-s+x,s+x)$ for all non-zero 
$y\in {\R}$ and the proposition is proved.
\end{proof}
\s
From the above construction it follows that the boundary
hypersurfaces $S(p)$ are strongly pseudoconvex.  Since then the 
smallest Stein domain containing a $S$-invariant neighborhood of 
$S(p)$ is $\Omega(p)\smallsetminus\diag (D_0)$,
the following is immediate.

\begin {cor}

For $p\in D_0\times D_0\smallsetminus \diag (D_0) $
every function $f$ which is holomorphic on some
neighborhood of the orbit $S\cdot p$ extends holomorphically
to $\Omega(p)\smallsetminus \diag (D_0)  $. An analogous statement 
is valid for $p\in D_\infty\times D_\infty\smallsetminus \diag (D_\infty) $.
\end {cor}
\s
Observe that the set $\rm bd_{gen}(D_0\times D_\infty)$ of generic 
boundary points, which was introduced in section 4.2, 
consists of the two $S$-orbits 
${\rm bd}(D_0)\times D_\infty\cup D_0\times {\rm bd}(D_\infty)$.
Let $z\in {\rm bd}(D_0)\times D_\infty$ 
(or $z\in D_0\times {\rm bd}(D_\infty)$ be such a boundary point.

\begin {cor}

Let $\widehat \Omega \subset Q_2\subset \C\P^1\!\times\C\P^1 $ be an
 $S$-invariant Stein domain  which
contains $D_0\times D_\infty$ and the boundary point $z$. 
Then $\widehat \Omega$ also  contains
$D_0\times \C\P^1\smallsetminus \diag (\C\P^1) $ 
(or $\C\P^1\times D_\infty\smallsetminus \diag (\C\P^1) $, respectively).

\end {cor}

\begin {proof}

Let $B$ be a ball around $z$ which is contained in 
$ \widehat \Omega $. For $p\in B(z)\cap D_\infty\times D_\infty$ 
sufficiently close to $z$ it follows that $S\cdot q 
\subset \widehat \Omega$ for all
$q\in B(z)\cap(D_\infty\!\times D_\infty)$.  
The result then follows from the previous Corollary.
\end {proof}
\s
If $\widehat \Omega $ is as in the above Corollary, then
the fibers of the projection of 
$\widehat \Omega\subset \C\P^1\!\times \C\P^1\to \P^1 $ can
be regarded as non-constant holomorphic curves
$f:\C\to \widehat \Omega $.  One says that a
complex manifold $X$ is Brody hyperbolic if there are
no such curves. 

\begin {cor} \label {nothyper}

If $\widehat \Omega$ is as above, then $\widehat \Omega$ 
is not Brody hyperbolic.
\end {cor} 
\s
A complex manifold $X$ is said to be Kobayashi hyperbolic whenever
the Kobayashi pseudo-metric is in fact a metric (see \cite {K}).
The pseudo-metric is defined in such a way that, if there
exists a non-constant holomorphic curve
$f:{\C}\to X$, then $X$ is not hyperbolic, i.e., Kobayashi
hyperbolicity is a stronger condition than Brody hyperbolic.
For an arbitrary semisimple group $G$ the domain $\OAG$ 
is indeed Kobayashi hyperbolic (\cite {H}, see $\S$5 for stronger
results). 
\s
The following is our main application of the existence of
$Q_2$-slices at generic points of $\text{\rm bd}(\OAG )$.

\begin {thm} \label {maximality}

A $G$-invariant, Stein and  Brody hyperbolic domain
$\widehat \Omega $ in $\GCKC$
which contains $\OAG $ is equal to $\OAG $.

\end {thm}

\begin {proof}

Arguing by contraposition, if $\OAG $ is strictly contained 
in a $G$-invariant Stein domain $\widehat \Omega $, then by Thm.~\ref{slice}
there exists a $Q_2$-slice at a generic boundary point 
$z \in \text{\rm bd}(\OAG )\cap \widehat \Omega$ 
with $Q_2\cap \widehat \Omega $ an $S$-invariant Stein domain 
properly containing the Akhiezer-Gindikin domain of $Q_2$. 
However, by Cor.~\ref{nothyper}, such a domain in $Q_2$ is not 
Brody hyperbolic.
\end {proof}
\b


\section {Hyperbolicity and the characterization of cycle domains} \label {V}

In this section  it is shown the Wolf cycle domains $\Omega _W(D)$
are Kobayashi hyperbolic. The above theorem then yields
their characterization (see \ref {description}).

\subsection {Families of hyperplanes} \label {V.1}

We start by proving a general result concerning
families of hyperplanes in projective space and their
intersections with locally closed subvarieties.
Since such a subvariety is usually regarded as being embedded
by sections of some line bundle, it is natural to regard
the projective space as the projectivization $\P(V^*)$ of
the dual space and a hyperplane in $\P(V^*)$ as a point in $\P(V)$.
\m
We will think of a subset $S\in \P(V)$ as parameterizing 
a family of hyperplanes in $\P(V^*)$.
A non-empty subset $S\subset \P(V)$ is said to have the normal
crossing property if for every $k\in \mathbb N$ there exist
$H_1,\ldots H_k\in S$ so that for every subset 
$I\subset \{ 1,\ldots ,k\} $  the intersection
$\bigcap _{i\in I}H_i$ is $|I| $-codimensional.  If 
$|I| \ge \dim_{\C}V$, this means that the intersection is empty.  
\m
In the sequel $\langle S\rangle$ denotes the complex linear span of $S$ in
$\P(V)$, i.e., the smallest plane in $\P(V)$ containing $S$.

\begin {prop}

A locally closed, irreducible real analytic subset $S$ with
$\langle S\rangle=\P(V)$ has the normal crossing property.

\end {prop}

\belowdisplayskip=2pt plus1pt minus2pt
\begin {proof}

We proceed by induction over $k$. For $k=1$ there is nothing to prove.
Given a set $\{ H_{s_1},\ldots ,H_{s_k}\} $ of hyperplanes
with the normal crossing property and a subset 
$I\subset \{s_1,\ldots ,s_k\} $, define
$$
\Delta _I:=\bigcap _{s\in I}H_{s}, 
\qquad {\cal H}(I):=\{ s\in S : H_s\supset \Delta _I \} 
\qquad {\Cl}_k:=\kern-1em\bigcup _{{\scriptstyle J\subset \{s_1,..,s_k\}
\atop
\scriptstyle \Delta _J\ne\varnothing}}\kern-1em{\cal H}(J).
$$
We wish to prove that $S\smallsetminus {\Cl}_k \ne\varnothing $.  For this,
note that each ${\cal H}(I)$ is a real analytic subvariety of $S$. Hence,
if $S={\Cl}_k$, then $S={\cal H}(J)$ for some $J$ with 
$\Delta _J \ne \varnothing $. However, 
$\{H\in \P(V^*):H\supset \Delta _J\} $ is a proper, linear plane ${\cal L}(J)$
of $\P(V)$.  Consequently, $S\subset {\cal L}(J)$, and this would 
contradict $\langle S\rangle=\P(V)$.  Therefore, there exists
$s \in S\smallsetminus {\Cl}_k$, or equivalently,
$\{ H_{s_1},\ldots ,H_{s_k},H_s\} $ has the normal crossing property.
\end {proof}
\belowdisplayskip=9pt plus1pt minus2pt
\s
It is known that if $H_1,\ldots, H_{2m+1}$ are hyperplanes 
having the normal crossing property, where $m=\dim_{\C}\P(V)$, then 
$\P(V^*)\smallsetminus \bigcup H_j$ is Kobayashi hyperbolic 
(\cite {D}, see also \cite {K} p.~137).  

\begin {cor}

If $S$ is a locally closed, irreducible and generating real analytic subset of 
$\P(V)$, then there exist hyperplanes $H_1,\ldots H_{2m+1}\in S$
so that the complement $\P(V^*)\smallsetminus \bigcup H_j$ is
Kobayashi hyperbolic.
\end {cor}
\s
Our main application of this result arises in the case where
$S$ is an orbit of the real form at hand.

\begin {cor}  \label {key}

Let $G^{\C}$ be a reductive complex Lie group, $G$ a real
form, $V^*$ an irreducible $G^{\C}$-representation space and
$S$ a $G$-orbit in $\P(V)$.  Then there exist
hyperplanes $H_1,\ldots ,H_{2m+1}\in S$ so that
$\P(V^*)\smallsetminus  \bigcup H_j$ is Kobayashi hyperbolic.

\end {cor}

\begin {proof}

From the irreducibility of the representation $V^*$, it follows
that $V$ is likewise irreducible and this, along with the identity
principle, implies that for $\langle S\rangle=\P(V)$.
\end {proof}

\subsection {Hyperbolic domains in $G^{\C}/K^{\C}$} \label {V.2}

As we have seen above, hypersurfaces $H$ in $\Omega =G^{\C}/K^{\C}$
which are invariant under the action of an Iwasawa-Borel group
$B$ play a key role in the study of $G$-invariant domains (see also 
\cite {H},\cite {HS},\cite {HW1}, and \cite {HW2}). In the sequel we
shall simply refer to such $H$ simply as a $B$-hypersurface.
\m
Recall that if $H_1,\ldots ,H_m$ are all of the  irreducible
$B$-hypersurfaces in $\GCKC$
and if $\bigcup _j(\bigcup _{g\in G}g(H_j))$ is removed from $\Omega $, then
the resulting set $\Omega _I$ is the Akhiezer-Gindikin domain
$\OAG$ (see Theorem \ref {equality}).  In particular, the resulting 
domain is non-empty. 
\s
Now if $H$ is just one (possibly not irreducible) 
$B$-hypersurface, then the 
set $\bigcup _{g\in G} g(H)=\bigcup _{g\in K} g(H) $ 
is closed and its complement in $\Omega= \GCKC$ is open. Let $\Omega _H$ be
the connected component of that open complement,
containing the chosen base point $x_0$.
It is likewise a non-empty $G$-invariant Stein domain in 
$\Omega= \GCKC$.  Here we shall prove that, if $G$ is not Hermitian,
any such $\Omega_H$ is Kobayashi hyperbolic. In the Hermitian
case one easily describes the situation where $\Omega _H$ is 
not hyperbolic. 
\m
Let $H$ be given as above and let $L$ be the line bundle which
it defines.  Let $\sigma_H$ be the corresponding 
section, i.e., $\{\sigma_H=0\}=H$. For convenience we may regard 
$L$ as an algebraic $G^{\C}$-bundle on a smooth, 
$G^{\C}$-equivariant projective compactification $X$ of $\GCKC$ 
to which the $B$-hypersurface $H$ extends.
\s
Note that  $\sigma_H $ is a $B$-eigenvector in $\Gamma (X,L)$.
Let $V_H\subset \Gamma (X ,L)$ be the irreducible $G^{\C}$-representation
space which contains $\sigma_H $. Define $\varphi _H:\Omega \to 
\P(V_H^*)$ to be the canonically associated $G^{\C}$-equivariant 
meromorphic map.  

\begin {lem}

The map 
$\varphi _H\vert_{\textstyle \Omega} :\Omega \rightharpoonup \P(V_H^*)$ 
is a regular morphism onto a quasi-projective $G^{\C}$-orbit 
$G^{\C}\cdot v_0^*=:\widetilde \Omega $.

\end {lem}

\begin {proof}

By definition $\varphi _H$ is $G^{\C}$-equivariant; in particular
its set $E$ of base points is $G^{\C}$-invariant.  Since $\Omega $ is
$G^{\C}$-homogeneous, $E=\varnothing $.
\end {proof}
\s
From now on we replace $\varphi _H$ by its restriction to 
$\Omega $ and only discuss that map. 
By definition every section $s\in V_H$ is the pull-back 
$\varphi _H^*(\widetilde s)$ of a hyperplane section. Thus,
there is a uniquely defined $B$-hypersurface $\widetilde H$
in $\P(V_H^*) $ with $\varphi _H^{-1}(\widetilde H)=H$.  Let
$\widetilde \Omega _{\widetilde H}\subset \P(V_H^*)$ be defined
analogously to $\Omega _H$, i.e.,  $\widetilde \Omega _{\widetilde H}= \P(V_H^*)\smallsetminus \bigcup _{g\in G} g(\widetilde H).$ 
Applying Cor.~\ref{key} to $\P(V_H^*)$ and 
$S:=G\cdot \widetilde H\subset \P(V_H)$, 
it follows that the domain
$\widetilde \Omega _{\widetilde H}$ is Kobayashi hyperbolic.
Further, the connected component of
$\varphi _H^{-1}(\widetilde \Omega _{\widetilde H})$ which contains the base point 
$x_0$ is just the original domain $\Omega _H$.
\m
If $\varphi $ has positive dimensional fibers, which indeed
can happen in the Hermitian case, then, since the connected
components of its fibers contain many holomorphic curves
$f:{\C}\to \Omega $, it follows that $\Omega _H$ is
not Kobayashi hyperbolic. 
\s
In the case of finite fibers, since
preimages under locally biholomorphic maps of hyperbolic manifolds
are hyperbolic,
the opposite is true.
\m

\begin {thm}

If the $\varphi _H$-fibers are finite, then $\Omega _H$ is
Kobayashi hyperbolic.

\end {thm}

\begin {cor} \label {kobayashi}

If $G$ is not of Hermitian type, then $\Omega _H$ is
Kobayashi hyperbolic.

\end {cor}

\begin {proof}

If $G$ is not of Hermitian type, then $K^{\C}$ is
dimension theoretically maximal in $G^{\C}$ and,
since $\varphi _H$ is non-constant, it follows that
it has finite fibers.
\end {proof}

\begin {thm} \label {cyclehyp}

The Wolf cycle domain $\Omega _W(D)$ of an open orbit $D$ of an
arbitrary real form $G$ of an arbitrary complex semisimple group
$G^\C $ in an arbitrary flag manifold $Z=G^\C/Q$ ist Stein
and Kobayashi hyperbolic.

\end {thm}

\begin {proof}

In (\cite {HW2}) it was shown that every Wolf cycle space
$\Omega _W(D)$ is the intersection of certain of the
$\Omega _H$. In the notation of (\cite {HW2}) such
an intersection is referred to as the associated Schubert
domain $\Omega _S(D)$.  Thus the cycle domains $\Omega _W(D)$ are Stein.
\s
If $G$ is not of Hermitian type, then, since it is contained 
in $\Omega _H$ for certain $B$-hypersurfaces $H$, Cor.~\ref {kobayashi}
implies that it is hyperbolic.
\m
If $G$ is of Hermitian type, then $\Omega _W(D)$ is either
the associated bounded symmetric domain ${\cal B}$, its complex
conjugate or, if $\Omega $ is non-compact, ${\cal B}\times \bar 
{\cal B}$ (\cite {W2},\cite {WZ1},\cite {WZ2},\cite {HW2}). Since
bounded domains are hyperbolic, this completes the proof.
\end {proof}  
\s
The following consequence of Cor.~\ref {nothyper} is from our
point of view the main technical result of this paper. 

\begin {thm} \label {notbrody}

If $\widehat \Omega $ is a $G$-invariant Stein domain in $
G^{\C}/K^{\C}$ which properly contains $\OAG$,
then $\widehat \Omega $ is not Brody hyperbolic.  

\end {thm}

\begin {proof}

Let $z\in \text{\rm bd}(\OAG)$ be a generic boundary
point in the sense of Prop. \ref {generic} which is also contained 
in $\widehat  \Omega $.
We may assume that $z$ is an optimal base point and that
$x_1=\exp(iA)x_0$ is in its closure as in Lemma \ref {optimal}.
Let $S^{\C}$ be as in Theorem \ref {slice} so that $Q_2=S^{\C}\cdot x_1$
is a $Q_2$-slice which contains $z$.  Since   
$Q_2\cap \widehat \Omega $ contains the Akhiezer-Gindikin domain of $Q_2$
($D_0\times D_\infty $ in the language of $\S 4.4$) as well as the 
boundary point $z$, it follows from Cor.~\ref{nothyper} that 
$\widehat \Omega $ is not Brody hyperbolic.
\end {proof}
\s
We now give a characterization of all Wolf cycle domains, including
the few exceptions mentioned above. For this recall that $D$ is
an open $G$-orbit in $Z=G^\C/Q$, $C_0$ the base cycle in $D$,
and $G^\C\cdot C_0=\Omega $ is the corresponding orbit in the
cycle space ${\cal C}^q(Z)$. 

\begin {thm} \label {description}
  
If $\Omega $ is compact, then
either $\Omega _W(D)$ consists of a single point
or $G$ is Hermitian and $\Omega _W(D)$ is either the
associated bounded symmetric domain ${\cal B}$ or
its complex conjugate $\bar {\cal B}$. If
$\Omega $ is non-compact, then, regarding $\Omega _W(D)$ as
a domain $G^\C/K^\C$, it follows that
$$\Omega _W(D)=\Omega _S(D)=\Omega _I=\OAG$$
for every open $G$-orbit in every $G^{\C}$-flag manifold
$Z=G^{\C}/Q$.

\end {thm}

\begin {proof}

The exceptional case where $\Omega $ is compact is discussed in
detail in the proof of Theorem~\ref{cyclehyp} and therefore
we restrict here to the non-compact case.
\s 
The statement $\Omega _W(D)=\Omega _S(D)$ is proved in (\cite {HW2}).
In $\S2$ above it is proved that $\OAG=\Omega _I$.  By definition
$\Omega _S(D)\supset \Omega _I=\OAG $.  
Since $\Omega _W(D)$ is 
Stein and hyperbolic (Theorem~\ref{cyclehyp}), by 
Theorem~\ref {notbrody} it follows that $\Omega _W(D)=\OAG$,
and all equalities are forced.
\end {proof}

\b
As a final remark we would like to mention (\cite {GM})
where cycle domains $C(\gamma )$ are introduced for every $G$-orbit
$\gamma \in Orb_Z(G)$. If $\kappa $ is the dual $K^\C$-orbit
to $\gamma $, then $C(\gamma )$ is the connected component
containing the identity of 
$\{ g\in G^\C: g(\kappa )\cap \gamma \ \text{\rm is compact
and non-empty} \}$.
Of course these spaces can be regarded in $G^\C/K^\C$ ,
and for open orbits $\gamma =D$ they are the same as the Wolf domains
$\Omega _W(D)$.
\m
Note that, by regarding the cycle domains as 
lying in the group $G^\C$ and going to the intersections,
one avoids the special considerations in the Hermitian case.
When the cycle spaces ${\cal B}$ and 
$\bar {\cal B}$ are regarded as subsets of ${\cal C}^q(Z)$,  
the intersection of
the associated spaces in $G^\C$ yields the domain
${\cal B}\times \bar {\cal B}$ in $\GCKC$.   
\m
Define
$$
C_Z(G):=\kern-.5em\bigcap _{\gamma \in {\rm Orb}_Z(G)}\kern-.5em  C(\gamma ).
$$ 
By Proposition 8.1 of \cite{GM},  
if
$x$ is in every cycle domain $\Omega _W(\widetilde D)$, where
$\widetilde D$ is an open $G$-orbit in $\widetilde  Z=G^\C/B$,
then $x\in C_Z(G)$ for every $Z$.
\m
The following result was proved in (\cite {GM}) for classical
and exceptional Hermitian groups via a case-by-case argument.

\begin {cor}

For $G$ an arbitrary real form of an arbitrary complex semisimple group 
$G^\C$ and $Z=G^\C/Q$ any $G^\C$-flag manifold,
it follows that $C_Z(G)=\OAG $.

\end {cor}

\begin {proof}

Using the above mentioned Prop. 8.1, the remark in
the Hermitian case and Theorem \ref {description}, 
it follows immediately that $\OAG \subset C_Z(G)$.
For the other inclusion, note that by definition
$C_Z(G)$ is contained in the intersection of the cycle
domains $\Omega _W(D)$ for $D$ an open $G$-orbit
in $Z$. By Theorem~\ref {description} this intersection
is again $\OAG $.
\end {proof}

\b
Fakult\" at f\" ur Mathematik \\
Ruhr--Universit\"at Bochum \\
D-44780 Bochum, Germany \\
{\tt gfels@cplx.ruhr-uni-bochum.de}\\
{\tt ahuck@cplx.ruhr-uni-bochum.de}
\end {document}